\documentclass[11pt]{article}
\usepackage[a4paper,bindingoffset=0.2in,%
            left=1in,right=1in,top=1in,bottom=1in,%
            footskip=.25in]{geometry}
\usepackage{amsmath,amsfonts,amssymb,latexsym,epsfig,subfig,float}
\usepackage{color,epsf}
\usepackage{booktabs,listings}
\usepackage{graphicx,stix}
\usepackage{alltt}

\newtheorem{remark}{Remark}

\newcommand{\psis}{{\mathcal S}}
\newcommand{\psib}{\chi}
\newcommand{\psiba}{\chi^*}
\renewcommand{\r}{\mathbb{R}}

\newcommand{\cM}{\r^D}
\newcommand{\hf}{F}
\newcommand{\cL}{{\mathcal L}}

\date{March 8, 2020}
\makeindex
\begin{document}
\title{Composition Methods for Dynamical Systems Separable into Three Parts}

\author{Fernando Casas\thanks{Email: \texttt{Fernando.Casas@uji.es}}
\and
Alejandro Escorihuela-Tom\`as\thanks{Email: \texttt{alescori@uji.es}}
}

\maketitle

\begin{abstract}

New families of fourth-order composition methods for the numerical integration of initial value problems defined by ordinary differential equations are proposed.
They are designed when the problem can be separated into three parts in such a way that each part is explicitly solvable. The methods
are obtained by applying different optimization criteria and preserve geometric properties of the continuous problem by construction. Different numerical examples
exhibit their improved performance with respect to previous splitting methods in the literature.

\vspace*{0.5cm}

\begin{center}
Institut de Matem\`atiques i Aplicacions de Castell\'o (IMAC) and  De\-par\-ta\-ment de
Ma\-te\-m\`a\-ti\-ques, Universitat Jaume I,
  E-12071 Cas\-te\-ll\'on, Spain.
\end{center}

\end{abstract}

\bigskip

\section{Introduction}
\label{sec.1}
Splitting methods are particularly useful for the numerical integration of ordinary differential equations (ODEs)
\begin{equation}   \label{ode.1}
   \dot{x} \equiv \frac{dx}{dt} = f(x),  \qquad x(t_0) = x_0 \in \mathbb{R}^D
 \end{equation}
when the vector field $f$ can be written as $f(x) =  \sum_{i=1}^n f_i(x)$, so that each subproblem 
\[
   \dot{x} = f_i(x), \qquad x(t_0) = x_0, \qquad i=1, \ldots, n
\]
is explicitly solvable, with solution $x(t) = \varphi_t^{[i]}(x_0)$. Then, by composing the different flows with appropriate chosen weights it is
possible to construct a numerical approximation to the exact solution $x(h)$ for a time-step $h$ of arbitrary order \cite{mclachlan02sm}. 
Although splitting methods have a long history in 
numerical mathematics and have been applied, sometimes with different names, in many different contexts (partial differential equations, quantum statistical mechanics, chemical physics, molecular dynamics, etc. \cite{glowinski16smi}), it is in the realm of
 {Geometric Numerical Integration} (GNI) where they play a key role, and in fact some of the most efficient geometric integrators are based on the
 related ideas of splitting and composition \cite{hairer06gni}. 
 
 In GNI the goal is
to construct numerical integrators in such a way that the approximations they furnish share one or several qualitative (often, geometric) properties with the
exact solution of the differential equation \cite{blanes16aci}. In doing so, the integrator has not only an improved qualitative behavior, but also allows for a significantly more 
accurate long-time integration than it is the case with general-purpose methods. In this sense, symplectic integration algorithms for Hamiltonian systems constitute
a paradigmatic example of geometric integrators \cite{sanz-serna94nhp,leimkuhler04shd}. Splitting and composition methods are widely used in GNI because the composition 
of symplectic (or volume preserving, orthogonal, etc.) transformations is again symplectic (volume preserving, orthogonal, etc., respectively).
 In composition methods the numerical scheme is constructed as the composition of several simpler integrators for the problem at hand, so as to 
 improve their accuracy.

When $f$ in (\ref{ode.1}) can be separated into two parts, very efficient splitting schemes have been designed and applied to solve a wide variety of problems 
arising in several fields, ranging from Hamiltonian Monte Carlo techniques to the evolution of the $N$-body gravitational problem in Celestial Mechanics (see \cite{hairer06gni,blanes16aci} and
references therein).

There are, however, relevant problems in applications where $f$ has to be decomposed into three or more parts in order to have subproblems that are explicitly solvable.
Examples include the disordered discrete nonlinear Schr\"odinger equation \cite{skokos14hot}, Vlasov--Maxwell equations in plasma 
physics~\cite{crouseilles15hsf}, the~motion of a charged particle
in an electromagnetic field according with the Lorentz force law~\cite{he15vpa} and problems in molecular dynamics~\cite{shang17anm}. 
In that case, although in principle methods of any order of accuracy can be built, the resulting
algorithms involve such a  large number of maps that they are not competitive in practice. It is the purpose of this paper to present an alternative class of efficient methods
for the problem at hand and compare their performance on some non-trivial physical examples than can be split into three parts.

The paper is structured as follows. We first review how splitting methods can be directly applied to get numerical solutions (Section \ref{sec.2}). Then the attention is turned
to the application of composition methods, and we get a family of 4th-order schemes obtained by applying a standard optimization procedure (Section \ref{sec.3}). In
Section \ref{sec.4} we show how standard splitting methods, when formulated as a composition scheme, lead to very competitive integrators, and also propose a different
optimization criterion for systems possessing invariant quantities. This allows us to get a new family of 4th-order schemes. All these integration algorithms are
subsequently tested in Section \ref{sec.5} on a pair of numerical examples. Finally, Section \ref{sec.6} contains some concluding remarks.

\section{First Approach: Splitting Methods}
\label{sec.2}

In what follows we assume that the vector field $f$ in (\ref{ode.1}) can be split into three parts, 
\begin{equation}  \label{3parts}
   f(x) = f_a(x) + f_b(x) + f_c(x) 
\end{equation}
in such a way that the exact $h$-flows $\varphi_h^{[a]}$, $\varphi_h^{[b]}$, $\varphi_h^{[c]}$, corresponding to $f_a$, $f_b$, $f_c$, 
respectively, can be computed exactly.

It is clear that the composition
\begin{equation}  \label{LT}
   \psib_h = \varphi^{[a]}_{h}\circ \, \varphi^{[b]}_{h}\circ \,
\varphi^{[c]}_{h} 
\end{equation}
(or any other permutation of the sub-flows) provides a first-order approximation to the exact solution  $x(h) = \varphi_h(x_0)$ of (\ref{ode.1}), i.e.,
\[
   \chi_h(x_0) =   \varphi_h(x_0) + \mathcal{O}(h^2),
\]
whereas the so-called Strang splitting
\begin{equation}  \label{strang}
  \mathcal{S}_h^{[2]} = \varphi^{[a]}_{h/2}\circ \, \varphi^{[b]}_{h/2}\circ \, 
\varphi^{[c]}_{h} \circ \, \varphi^{[b]}_{h/2} \circ \, \varphi^{[a]}_{h/2}
\end{equation}
leads to a second-order approximation. 

Higher order approximations to the exact solution of (\ref{ode.1}) can be obtained by generalizing (\ref{strang}), i.e., by
considering splitting schemes of the form
\begin{equation}   \label{split1}
  \psi_h^{[r]} = \varphi_{c_{s} h}^{[c]}  \circ \, \varphi_{b_{s} h}^{[b]} \circ \, \varphi_{a_{s} h}^{[a]} \circ  \cdots \circ \, 
    \varphi_{c_{1} h}^{[c]}  \circ \, \varphi_{b_{1} h}^{[b]} \circ \, \varphi_{a_{1} h}^{[a]},
\end{equation}
where the coefficients $a_i, b_i, c_i$, $i=1, \ldots, s$, are chosen to achieve a prescribed order of accuracy, say, $r$,
\begin{equation}  \label{oa1}
   \psi_h^{[r]}(x_0) =   \varphi_h(x_0) + \mathcal{O}(h^{r+1}) \quad \mbox{ as } \quad h \rightarrow 0.
\end{equation}
Requirement (\ref{oa1}) leads to a set of polynomial equations (the so-called 
{order conditions}), whose number and complexity grows enormously with the order. In particular, if $r=1$ (i.e., for a consistency method)
one has
\[
  \sum_{i=1}^{s}a_i=1, \qquad   \sum_{i=1}^{s}b_i=1, \qquad   \sum_{i=1}^{s}c_i=1.
\]
The specific number of order conditions is determined in fact by the dimension $c_k$ of the homogeneous subspace 
of grade $k$, $1 \le k \le r$, of the
free Lie algebra $\mathcal{L}(A,B,C)$ generated by the Lie derivatives $A, B, C$ corresponding to $f^{[a]}$, $f^{[b]}$, $f^{[c]}$, respectively
\cite{mclachlan02sm}. These dimensions are collected Table \ref{table.1} for $1 \le k \le 8$.

\begin{table}[H]
  \caption{{Number of order conditions to be satisfied by a splitting method of the form (\ref{split1}) at each order~$k$.}} \label{table.1}
  \centering
  \begin{tabular}{ccccccccc}
    \toprule
    Grade $k$ & 1 & 2 & 3 & 4 & 5 & 6 & 7  &  8 \\
    $c_k$  & 3 & 3 & 8 & 18 & 45 & 116 & 312 & 810  \\
    \bottomrule
  \end{tabular}
\end{table}

Thus, a splitting method (\ref{split1}) of order 4 requires solving $3+3+8+18 = 32$ order conditions and therefore the evaluation of at least a similar number
of sub-flows to have as many parameters as equations. This number can be reduced by considering time-symmetric methods, i.e., schemes verifying
\begin{equation}   \label{csym1}
  \psi_h^{[r]} \circ \, \psi_{-h}^{[r]} = \mathrm{id},
\end{equation}
where $\mathrm{id}$ is the identity map. Condition (\ref{csym1}) is verified by left-right palindromic compositions, i.e., if 
\[
   a_{s+1-i} = a_i, \quad b_{s+1-i} = b_i, \quad c_{s+1-i} = c_i,  \qquad i=1,2,\ldots
\]
in (\ref{split1}). Then all the conditions at even order are automatically satisfied. Thus, a symmetric method of order 4 requires solving $11$ order
conditions (instead of 32). Still, within this approach, one has to solve 56 polynomial equations to construct a method of order 6.

Methods of this class have been systematically analyzed in \cite{koseleff96eso}. In particular, it has been shown that 
if one aims to get schemes (\ref{split1}) of order 2 with the 
minimum number of maps, then the Strang splitting (\ref{strang}) is recovered. With respect to order 4, the following scheme was presented:
    \begin{equation}  \label{kosel1}
     \psi^{[4]}_{\tau}  =  \varphi^{[c]}_{c_1 \tau} \circ  \, \varphi^{[b]}_{b_1 \tau} \circ  \, \varphi^{[a]}_{a_1 \tau} \circ \, \varphi^{[b]}_{b_2 \tau} \circ \,
     \varphi^{[c]}_{c_2 \tau} \circ \, \varphi^{[b]}_{b_3 \tau} \circ \, \varphi^{[a]}_{a_2 \tau} \circ \, \varphi^{[b]}_{b_3 \tau} 
    \circ \, \varphi^{[c]}_{c_2 \tau} \circ \, \varphi^{[b]}_{b_2 \tau} \circ \, \varphi^{[a]}_{a_1 \tau} \circ \, \varphi^{[b]}_{b_1 \tau} \circ \, \varphi^{[c]}_{c_1 \tau} 
    \end{equation}
with
\[
  a_1 = w_1, \; a_2 = w_0, \; b_1 = b_2 =  \frac{w_1}{2}, \;  b_3 = \frac{w_0}{2}, \; c_1 = \frac{w_1}{2}, \; c_2 = \frac{w_0 + w_1}{2}
\]
and 
\[
  w_1 = \frac{1}{2 - 2^{1/3}}, \qquad w_0 = 1 - 2 w_1.
\]

In fact, 13 is the minimum number of maps required. More efficient schemes 
 involving 17 and 25 maps can also be found in \cite{koseleff96eso}.  For simplicity, we denote method (\ref{kosel1}) as
$ (c_1 b_1 a_1 b_2 c_2 b_3 a_2 b_3 c_2 b_2 a_1 b_1 c_1) $.
More~recently, in \cite{auzinger17psm} a method involving 21 maps of the form
\begin{equation}   \label{21maps}
 (a_1 b_1 c_1 a_2 b_2 c_2 a_3 b_3 c_3 a_4 b_4 a_4 c_3 b_3 a_3 c_2 b_2 a_2 c_1 b_1 a_1)
\end{equation}
has also been proposed and tested on several numerical examples.

\section{Second Approach: Composition Methods}
\label{sec.3}

As it is clear from the previous considerations, constructing high order splitting methods for systems separable into three parts requires solving
a large number of polynomial equations involving the coefficients, and this is a very challenging task in general.  For this reason, we turn our attention
to another strategy based on compositions of the first order $\psib_h = \varphi^{[a]}_{h}\circ \, \varphi^{[b]}_{h}\circ \,
\varphi^{[c]}_{h}$ and its  {adjoint},
\[
   \chi_h^* := (\chi_{-h})^{-1} =  \varphi^{[c]}_{h}\circ \, \varphi^{[b]}_{h} \circ \, \varphi^{[a]}_{h}
\]
with appropriately chosen weights. In other words, we look for integrators of the form
\begin{equation} \label{eq.2.1.1}
\psi_h = \psib_{\alpha_{2s} h}\circ
\psiba_{\alpha_{2s-1}h}\circ
\cdots\circ
\psib_{\alpha_{2}h}\circ
\psiba_{\alpha_{1}h}, \quad \mbox{ with } \quad (\alpha_1,\ldots,\alpha_{2s})
   \in \r^{2s}
\end{equation}
verifying in addition the time-symmetry condition $\alpha_{2s+1-i} = \alpha_i$ for all $i$. 

\begin{remark}
Methods of the form
\begin{equation} \label{eq.2.1.2}
  \psi_h = \psis_{\alpha_{m} h}^{[2]}\circ \,  \psis_{\alpha_{m-1} h}^{[2]} \circ
\cdots\circ
 \psis_{\alpha_{2} h}^{[2]}  \circ \, \psis_{\alpha_{1}h}^{[2]}  \quad \mbox{ with } \quad 
(\alpha_1,\ldots,\alpha_{m}) \in \r^{m}
\end{equation}
and $\alpha_{m+1-i} = \alpha_i$ (commonly referred in the literature as {symmetric compositions of symmetric methods}~\cite{mclachlan02sm}) verify a much
reduced number of order conditions and allows one to construct very efficient high-order schemes~\cite{hairer06gni}. Notice that, since 
the Strang splitting (\ref{strang}) verifies $\psis_h^{[2]} = \psib_{h/2} \circ \psiba_{h/2}$, then it is clear that when analyzing methods
(\ref{eq.2.1.1}) we also recover schemes of the form (\ref{eq.2.1.2}).
\end{remark}

\subsection{Analysis in Terms of Exponentials of Operators}

The analysis of the composition methods considered here can be conveniently done by considering the Lie operators associated with 
the vector fields involved and the graded free Lie algebra they generate.

As is well known, for 
each infinitely differentiable map
$g:\cM \longrightarrow \r$, the function $g(\varphi_h(x))$ admits an expansion of the form \cite{arnold89mmo,sanz-serna94nhp}
\[
g(\varphi_h(x)) = \exp(h \hf)[g](x) = g(x) + \sum_{k\geq 1}
\frac{h^k}{k!} \hf^k[g](x), \qquad x \in \cM,
\]
where $\hf$ is the Lie derivative associated with $f$,
\begin{equation}  \label{eq.1.1b}
  F  \equiv L_f = \sum_{i=1}^D \, f_i(x) \, \frac{\partial }{\partial x_i}.
\end{equation}

Analogously, for the basic method $\chi_h$ one can associate a series of linear operators so that \cite{blanes08sac}
\[
    g(\psib_h(x)) = \exp(Y(h))[g](x), \quad \mbox{ with }  \quad Y(h) = \sum_{k \ge 1} h^k Y_k
\]
for all functions $g$, whereas for its adjoint one has 
\[
   g(\psiba_h(x)) = \exp \big(-Y(-h) \big)[g](x).
\]   

Then the operator series associated with the integrator (\ref{eq.2.1.1}) is
\begin{eqnarray*}
\Psi(h) &=& \exp(-Y(-h\alpha_1)) \exp(Y(h \alpha_2)) \cdots
\exp(-Y(-h\alpha_{2s-1})) \exp(Y( h \alpha_{2s})).
\end{eqnarray*}

Notice that the order of the operators is the reverse of the maps in  (\ref{eq.2.1.1}) (\cite{hairer06gni} p. 88). Now, by
 repeated application of the Baker--Campbell--Hausdorff 
formula \cite{blanes16aci} we can express formally $\Psi(h)$ as the exponential of an operator $\tilde{F}(h)$,
\begin{equation}   \label{modH}
   \Psi(h) = \exp(\tilde{F}(h)), \qquad \mbox{ with } \qquad \tilde{F}(h) = \sum_{k \ge 1} h^k F_k,
\end{equation}   
$h^k F_k \in \cL_k$ for each $k\geq 1$
and
$\cL = \bigoplus_{k\geq 1} \cL_k$ is the graded free Lie algebra
generated by the operators  $\{ h Y_1,h^2 Y_2,h^3 Y_3,\ldots\}$, where, by consistency, $Y_1 = \hf$.
One has explicitly
\begin{eqnarray*}
 Y(h \alpha_i) & = & h \alpha_i Y_1 + (h \alpha_i)^2 Y_2 + (h \alpha_i)^3 Y_3 + \cdots \\
  -Y(-h \alpha_i) & = & h \alpha_i Y_1 - (h \alpha_i)^2 Y_2 + (h \alpha_i)^3 Y_3 - \cdots \\
\end{eqnarray*}
so that 
\begin{eqnarray}  \label{modF}
  \tilde{F}(h) & = &  h w_1 Y_1 + h^2 w_2 Y_2 + h^3 (w_3 Y_3 + w_{12} [Y_1, Y_2] ) \nonumber \\
   &  & + h^4 (w_4 Y_4 + w_{13} [Y_1, Y_3] + w_{112} [Y_1,[Y_1,Y_2]])  \\
   &  &  +  h^5 \big( w_5 Y_5 + w_{14} [Y_1,Y_4] + w_{113} [Y_1,Y_1,Y_3] \nonumber \\
   & & + w_{1112} [Y_1,Y_1,Y_1,Y_2] + w_{23} [Y_2,Y_3] + w_{212} [Y_2,Y_1,Y_2]  \big)+  \mathcal{O}(h^6),  \nonumber
\end{eqnarray}
where $[Y_2,Y_1,Y_2] \equiv [Y_2,[Y_1,Y_2]]$, etc, $[\cdot, \cdot]$ refers to the usual Lie bracket and 
$w_1, w_2, \ldots$ are polynomials in the coefficients $\alpha_i$. In particular, one has
\begin{equation}    \label{cond.or.4}
\begin{aligned}   
 &  w_1 = \sum_{i=1}^{2s} \alpha_i, \qquad\quad  w_2 = \sum_{i=1}^{2s} (-1)^i \alpha_i^2, \\
 &  w_3 = \sum_{i=1}^{2s} \alpha_i^3, \qquad\quad  w_4 = \sum_{i=1}^{2s} (-1)^i \alpha_i^4, \\
 & w_{12} = \frac{1}{2} \left(  \sum_{i=1}^{2s-1} (-1)^{i+1} \alpha_i^2 \sum_{j=i+1}^{2s} \alpha_j +  
        \sum_{i=1}^{2s-1}  \alpha_i \sum_{j=i+1}^{2s} (-1)^j \alpha_j^2 \right). 
\end{aligned}
\end{equation}

Thus, a time-symmetric 4th-order method has to satisfy only consistency ($w_1 = 1$) and the order conditions at order three,
$w_3 = w_{12} = 0$. Notice, then, that the minimum number of maps to be considered  is $s=3$. In that case the integrator reads
\begin{equation}  \label{s3}
  \psi_h = \chi_{\alpha_1} \circ \chi_{\alpha_2}^* \circ \chi_{\alpha_3} \circ \chi_{\alpha_3}^* \circ \chi_{\alpha_2} \circ \chi_{\alpha_1}^*
\end{equation}
and the unique (real) solution is given by
\[
  \alpha_1 = \alpha_2 = \frac{1}{2(2 - 2^{1/3})}, \qquad \alpha_3 = \frac{1}{2} - 2 \alpha_1.
\]

This scheme corresponds to the familiar triple-jump integrator \cite{yoshida90coh}
\begin{equation}  \label{triplejump}
 \psi_h = \psis_{\alpha h/2}^{[2]} \circ   \, \psis_{\beta h}^{[2]} \circ  \,  \psis_{\alpha h/2}^{[2]}  \qquad \mbox{ with } \quad \alpha = 1/(2- 2^{1/3}).
\end{equation}

If $\psib_h = \varphi^{[a]}_{h}\circ\varphi^{[b]}_{h}\circ
\varphi^{[c]}_{h}$, then $\psi_h$  involves 13 maps (the minimum number) and corresponds precisely to the splitting method (\ref{kosel1}). 

It is worth remarking that the order conditions (\ref{cond.or.4}) are general for {any} composition method of the form
(\ref{eq.2.1.1}), with independence of the particular basic first-order scheme $\chi_h$ considered, as long as $\chi_h$ and its
adjoint $\chi_h^*$ are included in the sequence. Thus, for instance, one might take the explicit Euler method as $\chi_h$ and the implicit
Euler method as $\chi_h^*$, {and also a symplectic semi-implicit method and its adjoint, 
leading to the symplectic partitioned Runge--Kutta schemes considered
in \cite{diele11esp}.}

\subsection{Composition Methods of Order 4}

Although one already gets a method of order 4 with only three stages, it is well known that the scheme (\ref{triplejump}) has large high-order
error terms. A standard practice to construct more efficient integrators consists in adding more stages in the composition and determine
the extra free parameters thus introduced according with some optimization criteria. Although assessing the quality of a given integration method
applied to all initial value problem is by no means obvious (the dominant error terms are not necessarily the same for different problems), several
strategies have been proposed along the years to fix these free parameters in the composition method (\ref{eq.2.1.1}).
Thus, in particular, one looks for solutions such that
the absolute value of the coefficients, i.e.,
\begin{equation}  \label{E1}
  E_1(\mbox{\boldmath $\alpha$}) =\sum_{i=1}^{2s} |\alpha_i| \qquad
\end{equation}
is as small as possible, the logic being that higher order terms in the expansion  (\ref{modF}) involve powers of these coefficients. In fact, methods with small values of 
$E_1(\mbox{\boldmath $\alpha$})$ usually have large stability domains and small error terms \cite{mclachlan02sm}.
In addition, for a number of problems, the dominant error term is precisely the coefficient $w_5$ multiplying $Y_5$ in the expansion (\ref{modF}), so
that it makes sense to minimize
\begin{equation}  \label{E2}
  E_{2}(\mbox{\boldmath $\alpha$})  =2s \, \big|\sum_{i=1}^{2s} \alpha_i^5 \big|^{1/4},
\end{equation}
for a given composition to take also into account the computational effort measured as the number $2s$ of basic schemes considered. 
Here, as in \cite{blanes06cmf}, we construct  symmetric methods  with small values of $E_1$ which, in addition, have also small values of $E_{2}$. For future reference, the corresponding values of the objective functions  for the triple-jump 
(\ref{triplejump}) are $E_1 = 4.40483$ and $E_2 = 4.55004$, respectively.

Next we collect the most efficient schemes we have obtained with $s=4,5,6$ by applying this strategy.

\paragraph{$s=4$ stages.}
The composition is 
\begin{equation}  \label{m.4}
  \psi_h = \chi_{\alpha_1} \circ \chi_{\alpha_2}^* \circ \chi_{\alpha_3} \circ \chi_{\alpha_4}^* \circ \chi_{\alpha_4} \circ  
  \chi_{\alpha_3}^* \circ \chi_{\alpha_2} \circ \chi_{\alpha_1}^*,
\end{equation}
and involves 17 maps when the basic scheme $\chi_h$ is given by (\ref{LT}). Now we have a free parameter, which~we take as $\alpha_1$. The  minima of both $E_1$ and $E_{2}$ are achieved at
approximately $\alpha_1 = 0.358$, and the resulting coefficients are collected in Table \ref{table.2} as method {XA}$_4$. In that case,
$E_1 = 2.9084$ and $E_{2} =  3.1527$.

\paragraph{$s=5$ stages.}
The resulting composition
\[
  \psi_h = \chi_{\alpha_1} \circ \chi_{\alpha_2}^* \circ \chi_{\alpha_3} \circ \chi_{\alpha_4}^* \circ \chi_{\alpha_5} \circ \chi_{\alpha_5}^* \circ
  \chi_{\alpha_4} \circ  
  \chi_{\alpha_3}^* \circ \chi_{\alpha_2} \circ \chi_{\alpha_1}^*
\]
involves 21 maps when applied to a system separable into three parts. Minimum values for $E_1$ and $E_2$ are achieved when
\[
    \alpha_1 = \alpha_2 = \alpha_3 = \alpha_4 = \frac{1}{2(4-4^{1/3})}, \qquad \alpha_5 = \frac{1}{2} - 4 \alpha_1.
\]

In consequence, the method can be written as
\[
 \psi_h = \psis_{\alpha h}^{[2]} \circ  \psis_{\alpha h}^{[2]} \circ   \psis_{\beta h}^{[2]} \circ    \psis_{\alpha h}^{[2]} \circ \psis_{\alpha h}^{[2]} 
\]
with  $\alpha = 2 \alpha_1$, $\beta = 2 \alpha_5$. Then  $E_1 = 2.3159$ and $E_{2} =     2.6111$. 
This method, denoted {XA}$_5$,  was first proposed in \cite{suzuki90fdo} and analyzed in detail in \cite{mclachlan02foh}. 

\paragraph{$s=6$ stages.}
Analogously we have considered a composition involving three free parameters (and 25 maps when $\chi_h$ is given by (\ref{LT})):
\begin{equation}  \label{m.6}
  \psi_h = \chi_{\alpha_1} \circ \chi_{\alpha_2}^* \circ \chi_{\alpha_3} \circ \chi_{\alpha_4}^* \circ \chi_{\alpha_5} \circ 
  \chi_{\alpha_6}^* \circ \chi_{\alpha_6} \circ
  \chi_{\alpha_5}^* \circ
  \chi_{\alpha_4} \circ  
  \chi_{\alpha_3}^* \circ \chi_{\alpha_2} \circ \chi_{\alpha_1}^*.
\end{equation}

 The proposed solution is collected in Table \ref{table.2} as method {XA}$_6$ leading to $E_1 = 2.0513$, $E_{2} =     2.4078$. Notice
how, by increasing the number of stages, it is possible to reduce the value of $E_1$ and $E_2$ as a measure of the efficiency of the schemes.
This integrator has been tested in the numerical integration of the so-called reduced $1 + 1/2$ Vlasov--Maxwell system \cite{bernier20smf}.

We could of course increase the number of stages. It turns out, however, that with $s=7$ one has the sufficient number of parameters to
satisfy all the order conditions up to order 6, resulting in a method of the form (\ref{eq.2.1.2}) \cite{yoshida90coh} involving 29 maps. More efficient 
6th-order schemes
can be obtained indeed by increasing the number of stages. Thus, in particular, with $s=9$ and $s=11$ one has the methods designed in
\cite{kahan97ccf} (37 maps) and \cite{sofroniou05dos} (45 maps), respectively, when the basic scheme is given by (\ref{LT}).

\begin{table}[H]
  \caption{\small{Fourth-order composition methods \textit{XA}$_s$  with $s$ stages minimizing $E_1$ and $E_2$. Method \textit{S}$_6$ 
    corresponds to the splitting method of (\cite{blanes02psp} Table 2) expressed as a composition scheme. \label{table.2}}}
  \centering
    \renewcommand\arraystretch{1.1}
    \begin{tabular}{ll}
    \toprule
      \multicolumn{2}{c}{\textit{XA}$_4$}\\
      \toprule
      $\alpha_1= 0.358$ &\qquad $\alpha_2= -0.47710242361717810834$  \\ 
      $\alpha_3= 0.35230499471528197958$ & \qquad $\alpha_4= 0.26679742890189612876$  \\ 
      \midrule
      \multicolumn{2}{c}{\textit{XA}$_5$}\\
      \midrule

      $\alpha_1= \alpha_2=\alpha_3=\alpha_4=      \displaystyle \frac{1}{2(4-4^{1/3})}$ & \qquad $\alpha_5= \frac{1}{2}-4\alpha_1$ \\ 
      \midrule
      \multicolumn{2}{c}{\textit{XA}$_6$}\\
      \midrule
      $\alpha_1= 0.16$ &\qquad $\alpha_2= 0.15$ \\
      $\alpha_3= 0.16$ & \qquad  $\alpha_4= -0.260672267225$ \\
      $\alpha_5=0.147945412322$ & \qquad $\alpha_6=0.142726854903$ \\
      \midrule
      \multicolumn{2}{c}{\textit{S}$_{6}$}\\
      \midrule
      $\alpha_1= 0.0792036964311957$ &\qquad $\alpha_2= 0.1303114101821663$ \\
      $\alpha_3= 0.22286149586760773$ & \qquad  $\alpha_4= -0.36671326904742574$ \\
      $\alpha_5=0.32464818868970624$ & \qquad $\alpha_6=0.10968847787674973$ \\
      \bottomrule
    \end{tabular}
\end{table}

\section{Third Approach: Splitting {via} Composition}
\label{sec.4}

We have already seen that
there exists a close relationship between composition methods of the form (\ref{eq.2.1.1}) and
splitting methods. This connection can be established more precisely as follows \cite{mclachlan95otn}.
Let us assume that $f$ in the ODE (\ref{ode.1}) can be split into {two} parts, $\dot{x} = f_a(x) + f_b(x)$, which each part explicitly solvable, and take 
$\chi_h =  \varphi^{[b]}_{h}\circ \, \varphi^{[a]}_{h}$. Then, the adjoint method reads $\chi_h^* = \varphi^{[a]}_{h} \circ \varphi^{[b]}_{h}$ and the composition 
(\ref{eq.2.1.1}) adopts the form
\begin{equation}
  \label{eq:MetAdj}
\psi_h = 
\big(\varphi^{[b]}_{\alpha_{2s} h}\circ \, 
\varphi^{[a]}_{\alpha_{2s}h}\big) \circ 
\big(\varphi^{[a]}_{\alpha_{2s-1} h}\circ \, 
\varphi^{[b]}_{\alpha_{2s-1}h}\big) \circ 
 \cdots\circ
\big(\varphi^{[b]}_{\alpha_{2} h}\circ \, 
\varphi^{[a]}_{\alpha_{2}h}\big) \circ 
\big(\varphi^{[a]}_{\alpha_{1} h}\circ  \,
\varphi^{[b]}_{\alpha_{1}h}\big) .
\end{equation}

Since $\varphi^{[i]}_{h}$, $i=a,b$ are {exact} flows, then they verify
$
 \varphi^{[i]}_{\beta h} \circ\varphi^{[i]}_{\delta h}= 
 \varphi^{[i]}_{(\beta+\delta) h}, 
$
and  (\ref{eq:MetAdj}) can be rewritten as the splitting scheme 
\begin{equation}
  \label{eq:splitting}
\psi_h = \varphi^{[b]}_{b_{s+1} h}\circ \, 
 \varphi^{[a]}_{a_{s}h}\circ \, \varphi^{[b]}_{b_{s} h}\circ \, 
 \cdots\circ \, 
 \varphi^{[b]}_{b_{2}h}\circ  \, 
 \varphi^{[a]}_{a_{1}h} \circ \, \varphi^{[b]}_{b_{1}h}
\end{equation}
if $b_{1}=\alpha_{1}$ and 
\begin{equation}
  \label{eq:abalpha}
a_j = \alpha_{2j} + \alpha_{2j-1}, \qquad\quad  
b_{j+1} = \alpha_{2j+1} + \alpha_{2j}, \qquad j=1,\ldots,s
\end{equation}
 (with $\alpha_{2s+1}=0$). Conversely, any integrator of the form (\ref{eq:splitting}) with 
 $\sum_{i=1}^{s} a_i = \sum_{i=1}^{s+1} b_i$ can be expressed in the form (\ref{eq.2.1.1}) 
 with $\chi_h =  \varphi^{[b]}_{h}\circ \, \varphi^{[a]}_{h}$ and 
\begin{equation}
  \label{eq:alphaab}
\begin{array}{l}
\alpha_{2s}=b_{s+1}, \\
\alpha_{2j-1}=a_{j}-\alpha_{2j}, \qquad
\alpha_{2j-2}=b_{j}-\alpha_{2j-1}, \qquad j=s,s-1,\ldots,1, \nonumber
\end{array}
\end{equation}
with $\alpha_{0}=0$ for consistency.
In consequence, {any} splitting method in principle designed for systems of the form $\dot{x} = f_a(x) + f_b(x)$ 
with no further restrictions on $f_a$ or $f_b$ can
be formulated as a composition~(\ref{eq.2.1.1})  which, in turn, can also be applied when $f$ is split into three (or more) pieces, $f = f_a + f_b + f_c$, by taking
$\chi_h =  \varphi^{[a]}_{h}\circ \, \varphi^{[b]}_{h} \circ  \, \varphi^{[c]}_{h}$. The performance will be in general different, since 
different optimization criteria are typically used. Notice that the situation is different, however, if~splitting methods of Runge--Kutta--Nystr\"om type are
considered.

A particularly efficient 4th-order splitting scheme designed for problems separated into two parts has been presented in (\cite{blanes02psp} Table 2) (method
\textit{S}$_{6}$) and will be used in our
numerical tests. It is a time-symmetric partitioned Runge--Kutta method of the form (\ref{eq:splitting}), since the role played by $f_a$ and
$f_b$ are interchangeable. When formulated as a composition method, it has six stages, i.e., it is of the form~(\ref{m.6}), with coefficients $\alpha_i$ listed in Table \ref{table.2}. For comparison, the corresponding values of $E_1$ and $E_2$ are
$E_1 = 2.4668$ and $E_{2} =     3.1648$.

\subsection*{An Optimization Criterion Based on the Error in Energy}

Very often, the class of problems to integrate are derived from a Hamiltonian function. In that case, Equation (\ref{ode.1}) is formulated as
\begin{equation}   \label{ham.1}
    \dot{q}_i = \frac{\partial H}{\partial p_i}, \qquad \dot{p}_i = -\frac{\partial H}{\partial q_i}, \qquad i = 1,\ldots, d
\end{equation}
so that $x = (q,p)^T$, $f = (\nabla_p H, -\nabla_q H)^T \equiv X_H$ and $H(q,p)$ is the Hamiltonian. The Lie derivative associated with $X_H$ verifies, for
any function $G: D \subset \mathbb{R}^{2d} \longrightarrow \mathbb{R}$,
\[
   L_{X_H} G = - \{ H, G \} = - \sum_{j=1}^{d} \left( \frac{\partial H}{\partial q_j}  \frac{\partial G}{\partial p_j}  - 
      \frac{\partial G}{\partial q_j}  \frac{\partial H}{\partial p_j} \right).
\]

In other words, $ \{ H, G \}$ is the Poisson bracket of $H$ and $G$. In this context, then, the Lie bracket of operators can be replaced by the real-valued Poisson
bracket of functions \cite{arnold89mmo}.

It is well known that the flow corresponding to (\ref{ham.1}) is symplectic and in addition preserves the total energy of the system. If $H$ can be
split as $H = A + B$, then $f^{[a]} = L_{X_A}$, $f^{[b]} = L_{X_B}$ and 
the splitting method (\ref{eq:splitting}) is also symplectic. Important as it is that the method shares this feature with the exact flow, 
one would like in addition that the energy be preserved as accurately as possible 
(since a numerical scheme cannot preserve both the symplectic form and the energy). A possible optimization criterion would be then to
select the free parameters in such a way that the error in the energy 
(or more in general, in the conserved quantities of the continuous system) is as small as possible. 

This criterion can be made more specific as follows \cite{blanes14nif}. 
First, we expand the modified Hamiltonian $\tilde{H}_h$ in the limit $h \rightarrow 0$ for a 4th-order splitting method
(\ref{eq:splitting}). A straightforward calculation shows that
\begin{equation}    \label{eq:tildeH}
  \begin{aligned}
   & \tilde{H}_h  =   H  +h^4 k_{5,1} \{A,A,A,A,B\}   +h^4 k_{5,2} \{B,A,A,A,B\}  + h^4 k_{5,3} \{A,A,B,A,B\}   \\
  &  + h^4 k_{5,4} \{A,B,B,A,B\}  + h^4 k_{5,5} \{B,A,B,A,B\} + h^4 k_{5,6} \{B,B,B,A,B\}+ \\
  & + h^5 \sum_{j=1}^9 k_{6,j} E_{6,j} + \mathcal{O}(h^6),
 \end{aligned}  
\end{equation}
where $k_{i,j}$ are polynomials in the coefficients $a_j$, $b_j$, $\{A,A,A,A,B\}$ refers to the iterated Poisson bracket
$ \{A,\{A,\{A,\{A,B\} \} \} \}$, and $E_{6,j}$ are (independent) Poisson brackets involving 6 functions $A$ and $B$.

Now the Lie formalism allows one to get the Taylor expansion of the energy after one time-step  (\cite{sanz-serna94nhp}, Section 12.2) as
\[
H(q_{i+1},p_{i+1}) = \exp(-h{\cal L}_{\tilde H_h})H(q_i,p_i)  = H(q_i,p_i) - h {\cal L}_{\tilde H_h}H(q_i,p_i)
+ \frac{1}{2}h^2 {\cal L}_{\tilde H_h}^2H(q_i,p_i)+\cdots,
\]
where ${\cal L}_{\tilde H_h}(\cdot) = \{\tilde H_h,\cdot\}$.

An elementary calculation shows that 
\[
\begin{aligned}
  & H(q_{i+1},p_{i+1}) - H(q_i,p_i) = h^5 \big( k_{51} E_{61} + (k_{51}-k_{53}) E_{62} +  (k_{52}-k_{53}) E_{63}  + k_{54} E_{64} \\
  & \;\; + (k_{52} - \frac{1}{3} k_{53}) E_{65} + (k_{55} - \frac{1}{3} k_{53}) E_{66} + (k_{55} + k_{54}) E_{67} + (k_{56} - k_{54}) E_{68}  + k_{56} E_{69} \big) \\
  & \;\; + \mathcal{O}(h^6). 
\end{aligned}
\]

Thus, for small $h$, 
\begin{equation}  \label{delta.1}
\begin{aligned}
 &  \Delta \equiv  k_{51}^2 + (k_{51}-k_{53})^2 + (k_{52}-k_{53})^2 +  k_{54}^2 + (k_{52} - \frac{1}{3} k_{53})^2 \\
 & \quad + (k_{55} - \frac{1}{3} k_{53})^2 + (k_{55} + k_{54})^2 + (k_{56} - k_{54})^2 +  k_{56}^2
\end{aligned} 
\end{equation}
can be taken as a measure of the energy error, and consequently, 
\begin{equation}  \label{E3}
  E_{3} =2s \,  \Delta^{1/4}
\end{equation}
constitutes a possible objective function to minimize. The previous analysis can be also carried out for a composition method (\ref{eq.2.1.1}), resulting in
\begin{equation}  \label{delta.2}
  \Delta = w_{5}^2 + w_{14}^2 + w_{113}^2 + w_{1112}^2 + w_{23}^2 + w_{212}^2.
\end{equation}  

The $s$-stage methods \textit{XB}$_{s}$ whose coefficients are collected in Table \ref{table.3} have been obtained by minimizing $E_3$ with (\ref{delta.2}) and
in addition provide small values for (\ref{delta.1}) when applied with $\chi_h = \varphi_h^{[b]} \circ \varphi_h^{[a]}$.

We should emphasize again that, although methods \textit{XB}$_{s}$ have been obtained by minimizing~(\ref{delta.2}), and thus the local
error in the energy, their applicability is by no means limited to Hamiltonian systems. As a matter of fact, both classes of schemes 
\textit{XA}$_{s}$ and \textit{XB}$_{s}$ can be used with \textit{any}  first-order basic method and its adjoint. Their efficiency may depend, of course,
of the type of problem one is approximating and the particular basic scheme taken to form the composition. Moreover, due to the close
relationship between symplectic and composition methods, these schemes can also be seen as 
symplectic partitioned Runge--Kutta methods that, in contrast to splitting schemes, do not require the knowledge of the solution of the elementary flows.

\begin{table}[H]
  \caption{\small{Fourth-order composition methods \textit{XB}$_{s}$ with $s$ stages minimizing $E_3$. \label{table.3}}}
  \centering
  \renewcommand\arraystretch{1.1}
  \begin{tabular}{ll}
  \toprule
    \multicolumn{2}{c}{\textit{XB}$_{4}$}\\
    \midrule
    $\alpha_1= 0.1728230091082606$ & \qquad $\alpha_2= 0.43074941762060376$ \\
    $\alpha_3= -0.5742238363039501$ & \qquad $\alpha_4= 0.4706514095750858$  \\ 
    \midrule
    \multicolumn{2}{c}{\textit{XB}$_{5}$}\\
    \toprule
    $\alpha_1= 0.08967664078837478$ &\qquad $\alpha_2= 0.16032335921162522$ \\
    $\alpha_3= 0.29632291754168816$ & \qquad $\alpha_4= -0.49421908717228863$ \\
    $\alpha_5=0.44789616963060047$ & \\ 
    \midrule
    \multicolumn{2}{c}{\textit{XB}$_{6}$}\\
    \midrule
    $\displaystyle \alpha_1= \frac{1}{20}$ &\qquad $\displaystyle \alpha_2= \frac{71}{660}$ \\[0.3cm]
    $\displaystyle \alpha_3= \frac{47}{330}$ & \qquad $\displaystyle \alpha_4= \frac{37}{165}$ \\[0.3cm]
    $\displaystyle \alpha_5=-\frac{313}{660}$ & \qquad $\displaystyle \alpha_6=\frac{5}{11}$ \\
    \bottomrule
  \end{tabular}
\end{table}

\section{Numerical Examples}
\label{sec.5}

Although optimization criteria based on the objective functions $E_1$, $E_2$ and $E_3$ allow one in principle to construct  efficient 
composition schemes, it is clear that their overall performance depends very much on the particular problem considered, the initial conditions, etc. 
It is, then, worth considering some illustrative numerical examples to test the methods proposed here with respect to other integrators 
previously available in the literature. In particular, we take as representatives the splitting method (\ref{21maps}) designed in \cite{auzinger17psm}
for problems separated into three parts (referred to as \textit{ABC}$_{21}$ in the sequel) and the splitting scheme of \cite{blanes02psp} considered as
a composition (\ref{eq.2.1.1}) (referred as \textit{S}$_6$ in Table \ref{table.2}).

When a specific composition method (\ref{eq.2.1.1}) is applied to a particular problem of the form $\dot{x} = f_a + f_b + f_c$
and the first-order method is $\psib_h = \varphi^{[a]}_{h}\circ \, \varphi^{[b]}_{h}\circ \, \varphi^{[c]}_{h}$, the implementation is in fact very similar as for
a splitting method of the form (\ref{21maps}). Thus, in particular, for the integrator (\ref{m.6}) one has to apply  the following
procedure for the time step $x_n \longmapsto x_{n+1}$,
where one has to take into account the symmetry of the coefficients: $\alpha_{12} = \alpha_1$, etc. {and $s=6$}:
\[
\begin{array}{l}\label{scheme.1}
  y = x_n    \\
  \mbox{do }   j=1:6  \\
     \quad y = \varphi_{\alpha_{2j-1} h}^{[a]} y \\
     \quad y = \varphi_{\alpha_{2j-1} h}^{[b]} y \\
     \quad  \widetilde{\alpha} = \alpha_{2j-1} + \alpha_{2j} \\
     \quad y = \varphi_{\widetilde{\alpha}h}^{[c]} y \\
     \quad y = \varphi_{\alpha_{2j} h}^{[b]} y \\
     \quad y = \varphi_{\alpha_{2j} h}^{[a]} y \\
 \mbox{end } \\
 x_{n+1} = y
\end{array}        
\]
It is worth remarking that the examples considered here have been chosen because they admit an straightforward separation into three
parts that are explicitly solvable and thus may be used as a kind of testing bench to illustrate the main features of the proposed algorithms. Of course,
many other systems could also be considered, including non linear oscillators and the time integration of Vlasov-Maxwell equations 
\cite{crouseilles15hsf,bernier20smf}. In addition, the general technique proposed in \cite{tao16esa} for obtaining explicit symplectic approximations of
non-separable Hamiltonians provides in a natural manner examples of systems separable into three parts.

\subsection{Motion of a Charged Particle under Lorentz Force}

Neglecting relativistic effects, the evolution of a  particle of mass $m$ and charge $q$ in a given electromagnetic field is described by
the Lorentz force as
\begin{equation}   \label{lorentz.1}
   m \, \ddot{\mathbf{x}} = q \, (\mathbf{E} + \dot{\mathbf{x}} \times \mathbf{B}),
\end{equation}   
where $\mathbf{E}$ and $\mathbf{B}$ denote the electric and magnetic field, respectively. In terms of position and velocity, the equation of motion
(\ref{lorentz.1}) can be restated as
\begin{equation}   \label{lorentz.2}
 \begin{aligned}
   &  \dot{\mathbf{x}} = \mathbf{v} \\
   &  \dot{\mathbf{v}} = \frac{q}{m}  \mathbf{E} + \omega \, \mathbf{b} \times \mathbf{v}
  \end{aligned}
\end{equation}
where $\omega = -q B/m$ is the local cyclotron frequency, $B = \|\mathbf{B}\|$ and $\mathbf{b} = \mathbf{B}/B$ is the unit vector in the direction of
the magnetic field. For simplicity, we assume that both $ \mathbf{E}$ and $ \mathbf{B}$ only depend on the position $\mathbf{x}$. 

System (\ref{lorentz.2}) can be split into three parts in such a way that (a) each subpart is explicitly solvable and (b) the volume form in the space
$(\mathbf{x}, \mathbf{v})$ is exactly preserved \cite{he15vpa,he16hov}:
\begin{eqnarray}   \label{lorentz.3}
  \frac{d}{dt} \left(  \begin{array}{c}
  				\mathbf{x}  \\
				\mathbf{v} 
			  \end{array}  \right)  & = &  
		 \left(  \begin{array}{c}
  				\mathbf{v}  \\
				0 
			  \end{array}  \right)	  + 
		 \left(  \begin{array}{c}
  				0  \\
				\frac{q}{m} \mathbf{E}(\mathbf{x}) 
			  \end{array}  \right)	  	+ 
		 \left(  \begin{array}{c}
  				0  \\
				\omega(\mathbf{x})  \mathbf{b}(\mathbf{x}) \times \mathbf{v}
			  \end{array}  \right)	  \nonumber 	\\
			 & = & f^{[a]}(\mathbf{x}, \mathbf{v}) +  f^{[b]}(\mathbf{x}, \mathbf{v}) + f^{[c]}(\mathbf{x}, \mathbf{v}). 
\end{eqnarray}

The corresponding flows with initial condition $(\mathbf{x_0}, \mathbf{v}_0)$ are given by
\begin{eqnarray}   \label{lorentz.4}
   &  &  \varphi_t^{[a]}: \left\{  \begin{array}{l}
   					\mathbf{x}(t) = \mathbf{x}_0 + t \, \mathbf{v}_0 \\
					\mathbf{v}(t) = \mathbf{v}_0
				\end{array}  \right., \qquad \qquad	
   \varphi_t^{[b]}: \left\{  \begin{array}{l}
   					\mathbf{x}(t) = \mathbf{x}_0  \\
					\mathbf{v}(t) = \mathbf{v}_0 + t \, \frac{q}{m} \, \mathbf{E}(\mathbf{x}_0) 
				\end{array}  \right.  \nonumber \\	
  &  &  \varphi_t^{[c]}: \left\{  \begin{array}{l}
   					\mathbf{x}(t) = \mathbf{x}_0  \\
					\mathbf{v}(t) = \exp(t \omega(\mathbf{x}_0) \hat{\mathbf{b}}_0)  \mathbf{v}_0
				\end{array}  \right.  
\end{eqnarray}
where $\hat{\mathbf{b}}_0 \equiv \hat{\mathbf{b}}(\mathbf{x}_0)$ is the skew-symmetric matrix
\[
   \hat{\mathbf{b}}(\mathbf{x}) = \left(  \begin{array}{ccc}
   				0  &  -b_3(\mathbf{x})  &  b_2(\mathbf{x}) \\
			  b_3(\mathbf{x})   &  0  &  -b_1(\mathbf{x}) 	\\
			  -b_2(\mathbf{x})  &  b_1(\mathbf{x}) &  0
			      \end{array}  \right)
\]
associated with $\mathbf{b}(\mathbf{x}) = (b_1(\mathbf{x}), b_2(\mathbf{x}), b_3(\mathbf{x}))^T$.   

As in \cite{he15vpa}, we consider a static, non-uniform electromagnetic field
\begin{equation}   \label{loretnz.5}
   \mathbf{E} = - \nabla V = \frac{0.01}{r^3} (x \, \mathbf{e}_x + y \ \mathbf{e}_y), \qquad   \mathbf{B} = \nabla \times \mathbf{A} = r \, \mathbf{e}_z
\end{equation}
derived from the potentials
\[
   V = \frac{0.01}{r}, \qquad \mathbf{A} = \frac{r^2}{3} \, \mathbf{e}_{\theta}
\]
respectively, in cylindrical coordinates $(r, \theta, z)$ and with the appropriate normalization. Then, it can be shown that both the angular momentum
and energy
\[
   L = r^2 \dot{\theta} + \frac{r^3}{3}, \qquad \qquad 
    H = \frac{1}{2} \|\mathbf{v}\|^2 + \frac{0.01}{r}
\]
are invariants of the problem \cite{he15vpa}. 
   
\begin{figure}[!h]
  \begin{center}
    \subfloat[]{
      \includegraphics[width=7.5cm]{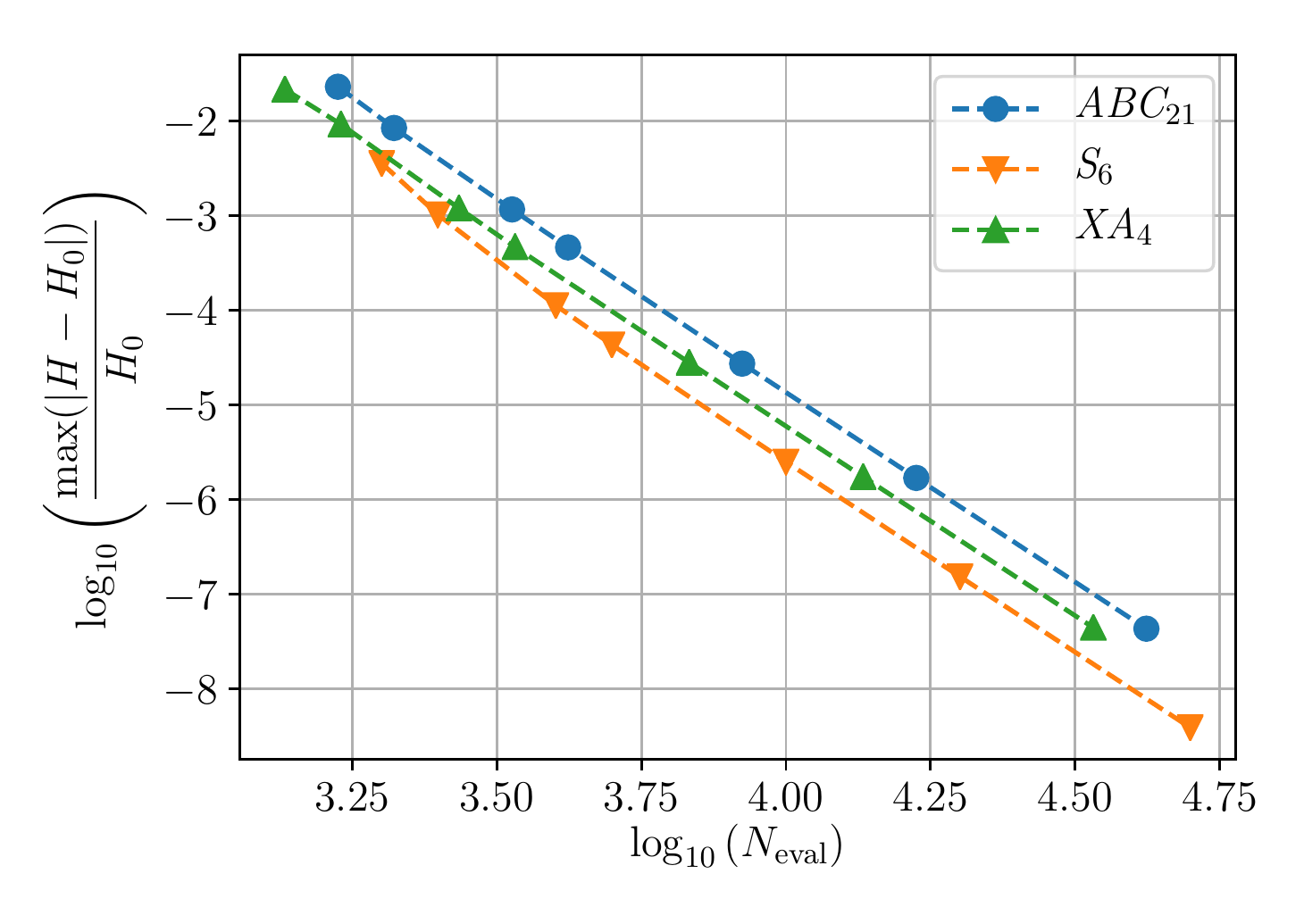}}
    \subfloat[]{
      \includegraphics[width=7.5cm]{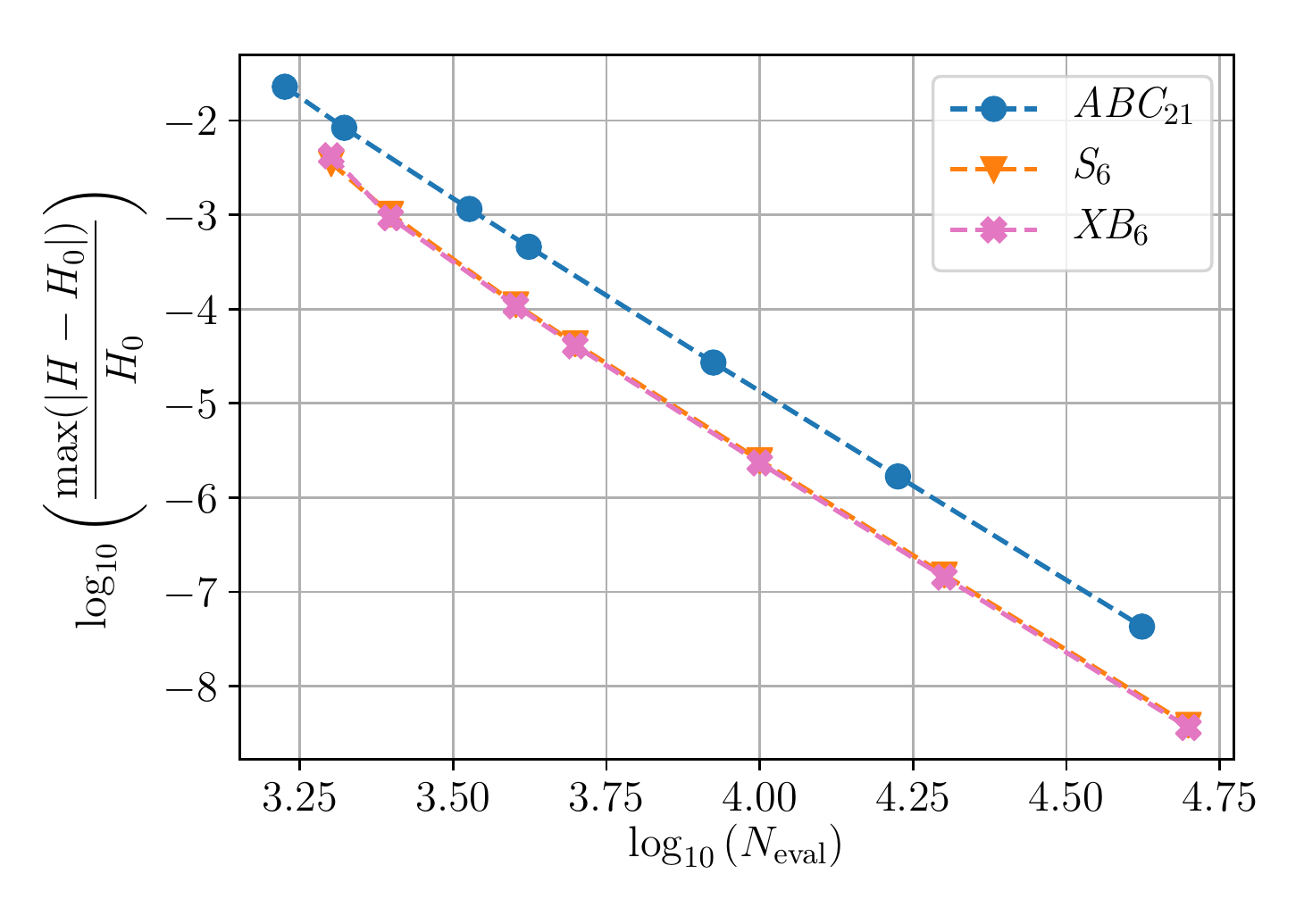}}
    \\ 
    \subfloat[]{
      \includegraphics[width=7.5cm]{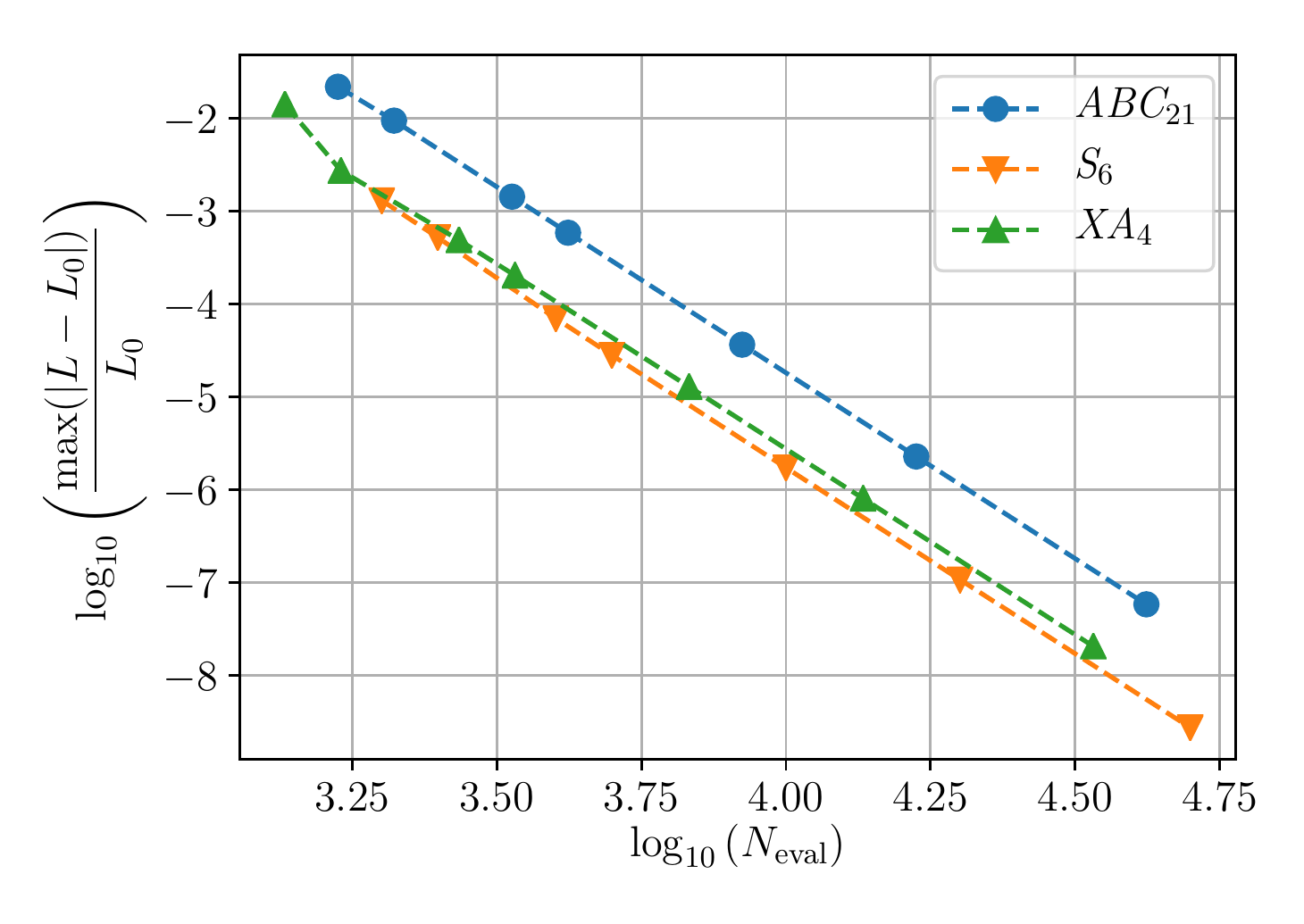}}
    \subfloat[]{
    \includegraphics[width=7.5cm]{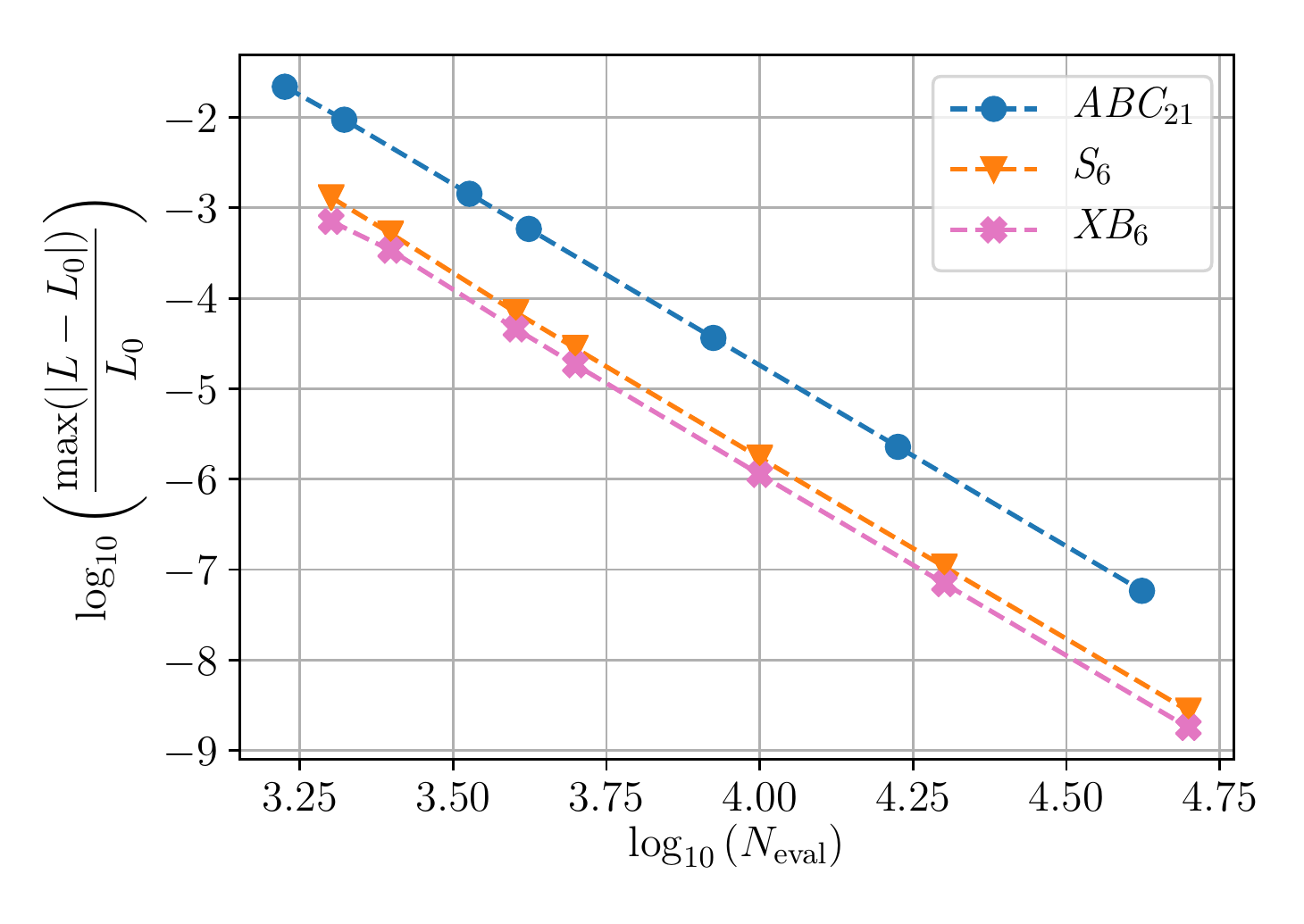}}
\caption{\small Relative error in conserved quantities due to each of the best numerical methods tested for charged particle under Lorentz force. (\textbf{a}) Relative error in energy for \textit{XA}$_4$ compared to \textit{ABC}$_{21}$ and  \textit{S}$_6$. (\textbf{b}) Relative error in energy for \textit{XB}$_6$ compared to \textit{ABC}$_{21}$ and  \textit{S}$_6$. (\textbf{c})Relative error in angular momentum for \textit{XA}$_4$ compared to \textit{ABC}$_{21}$ and  \textit{S}$_6$. (\textbf{d}) Relative error in angular momentum for \textit{XB}$_6$ compared to \textit{ABC}$_{21}$ and  \textit{S}$_6$.
\label{fig:EM1}}
  \end{center}
\end{figure}

With $q=-1$, $m=1$ and starting from the initial position $\mathbf{x}_0 = ( 0, -1, 0)^T$ with initial velocity $\mathbf{v}_0 = (0.10, 0.01, 0)$, 
we integrate with the different numerical
schemes until the final time $t_f = 200$ and compute the error in energy and angular momentum along the integration interval. As reference solution we take
the output generated by the standard routine DOP853 based on a Runge--Kutta method  of order 8 with local error estimation and step size control
(with a very stringent tolerance) \cite{hairer93sod}. In this way, we obtain Figure \ref{fig:EM1} (top and bottom, respectively), where this error is depicted
in terms of the number of the computed sub-flows (by taking different time-steps). For clarity, here and in the
sequel, in the left panel we include the results attained by the most efficient \textit{XA}$_s$  method, whereas the right panel corresponds to the 
\textit{XB}$_s$ schemes. For reference and comparison, we include in all cases the splitting method (\ref{21maps}) proposed in \cite{auzinger17psm} (denoted here
as \textit{ABC}$_{21}$)
and the scheme \textit{S}$_6$, whose coefficients are collected in Table \ref{table.2}.
   
We notice that applying the composition methods proposed here leads to more accurate results than the direct approach based on the splitting methods of Section
\ref{sec.2} with the same computational cost, and that the new scheme \textit{XB}$_6$ is slightly more efficient that the the splitting scheme \textit{S}$_6$ (the
remaining composition methods of Tables \ref{table.2} and \ref{table.3} provide results between \textit{ABC}$_{21}$ and the best composition depicted here).

\begin{figure}[!h]
  \begin{center}
    \subfloat[]{
      \includegraphics[width=7.5cm]{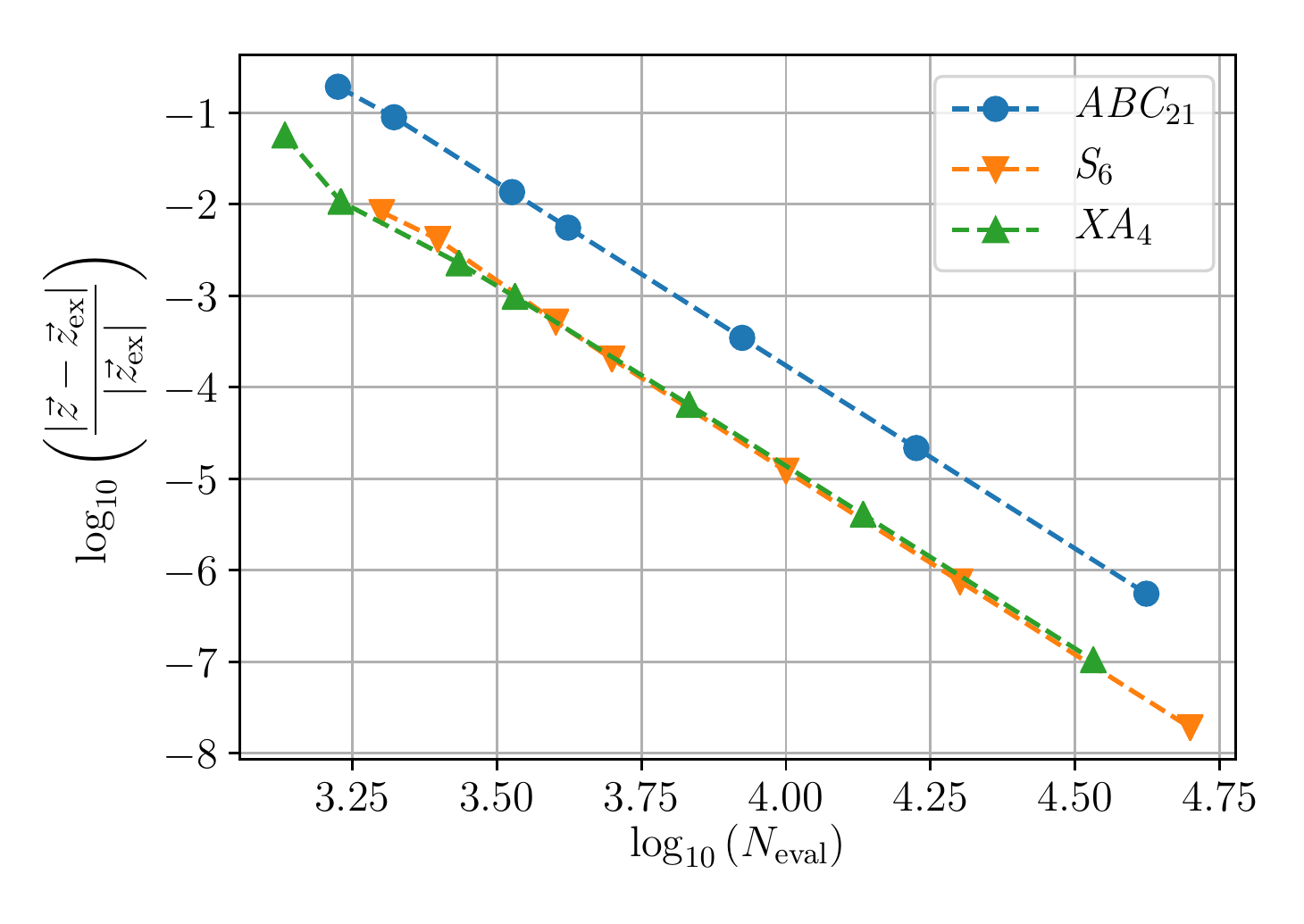}}
    \subfloat[]{
      \includegraphics[width=7.5cm]{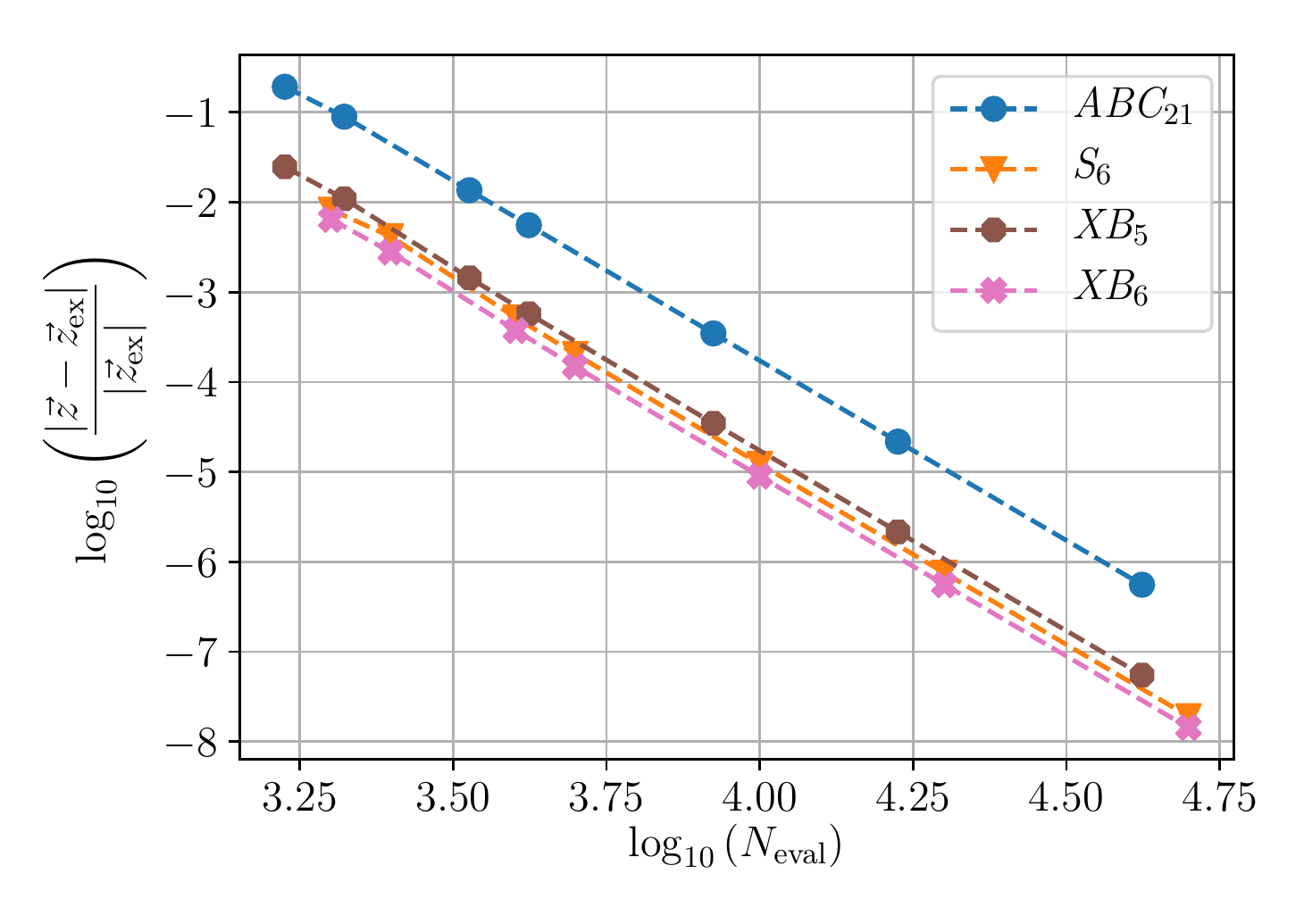}}   
\caption{\small Relative error in the $(\mathbf{x},\mathbf{v})$ space for charged particle under Lorentz force. The notation is the same as in Figure \ref{fig:EM1}. (\textbf{a}) Relative error in the $(\mathbf{x},\mathbf{v})$ space for \textit{XA}$_4$ compared to \textit{ABC}$_{21}$ and  \textit{S}$_6$. (\textbf{b}) Relative error in the $(\mathbf{x},\mathbf{v})$ space for \textit{XB}$_5$ and \textit{XB}$_6$ compared to \textit{ABC}$_{21}$ and  \textit{S}$_6$.
  \label{fig:EM3}}
\end{center}
\end{figure}

In Figure \ref{fig:EM3} we show the corresponding results obtained by each method for the error in the  $(\mathbf{x},\mathbf{v})$~space.

One should notice that, although this system is Hamiltonian, the Hamiltonian function is not separable into kinetic plus potential energy, and thus general
symplectic Runge--Kutta methods cannot be explicit \cite{he16hov}. In order to use explicit methods, one has to split the system into three parts. 
On the other hand, all the methods
tested here are volume-preserving in the  $(\mathbf{x},\mathbf{v})$ space, just as the exact flow.

\subsection{Disordered Discrete Nonlinear Schr\"odinger Equation}

The Hamiltonian of the disordered discrete nonlinear Schr\"odinger equation (DDNLS) 
\begin{equation}  \label{dnls.1}
  \mathcal{H} = \sum_{j} \left( \epsilon_j |\psi_j|^2 + \frac{\beta}{2} |\psi_j|^4 - (\psi_{j+1} \overline{\psi}_j + \overline{\psi}_{j+1} \psi_j) \right)
\end{equation}
describes a one-dimensional chain of couples nonlinear oscillators \cite{skokos14hot}. Here the sum extends over $N$ oscillators, $\psi_j$ are complex variables, 
$\beta \ge 0$ stands for the nonlinearity strength and the random energies $\epsilon_j$ are chosen uniformly from the interval $[-W/2, W/2]$,
where $W$ is related with the disorder strength. This model has two invariants: the energy (\ref{dnls.1}) and the norm
\[
   S = \sum_j |\psi_j|^2,
\]
and has been used to determine how the energy spreads in disordered systems \cite{kopidakis08aow}. Rather than analyzing the rich dynamics this 
system possesses, our interest here is to use (\ref{dnls.1}) as a non-trivial test bench for the integrators we presented in previous sections. 
By introducing the new (real) generalized coordinates and momenta $(q_j, p_j)$ related with $\psi_j$ through
\[
    \psi_j = \frac{1}{\sqrt{2}} (q_j + i p_j), \qquad  \overline{\psi}_j = \frac{1}{\sqrt{2}} (q_j - i p_j),
\]
the Hamiltonian function (\ref{dnls.1}) can be written as
\begin{equation}   \label{dnls.2}
  H = \sum_{j=1}^N \left( \frac{\epsilon_j}{2} (q_j^2 + p_j^2) + \frac{\beta}{8} (q_j^2 + p_j^2)^2 -  p_{j+1} \, p_j - q_{j+1} \, q_j \right)
\end{equation}  
in such a way that is the sum of three explicitly solvable parts, $H = A + B + C$, with
\[
  A = \sum_{j=1}^N \left( \frac{\epsilon_j}{2} (q_j^2 + p_j^2) + \frac{\beta}{8} (q_j^2 + p_j^2)^2 \right), \quad 
  B = - \sum_{j=1}^N   p_{j+1} \, p_j, \quad 
  C = -  \sum_{j=1}^N q_{j+1} \, q_j. 
\]  

The corresponding flows are given, respectively, by
\begin{eqnarray}   \label{dnls.3}
     \varphi_t^{[a]}: &  & \left\{  \begin{array}{l}
   					q_j(t) = q_j(t_0) \cos(a_j t) + p_j(t_0) \sin(a_j t)  \\
					p_j(t) = -q_j(t_0) \sin(a_j t) + p_j(t_0) \cos(a_j t) 
				\end{array}  \right., \nonumber 	\\
   \varphi_t^{[b]}: &  &  \left\{  \begin{array}{l}
   					q_j(t) = q_j(t_0) - t (p_{j-1}(t_0) + p_{j+1}(t_0))  \\
					p_j(t) = p_{j}(t_0) 
				\end{array}  \right.   \\	
   \varphi_t^{[c]}: &  &  \left\{  \begin{array}{l}
   					q_j(t) = q_j(t_0)  \\
					p_j(t) = p_j(t_0) + t (q_{j-1}(t_0) + q_{j+1}(t_0)) 
				\end{array}  \right.  \nonumber
\end{eqnarray}
with $a_j = \epsilon_j + \beta (q_j^2 + p_j^2)/2$.

\begin{figure}[!h]
  \begin{center}
    \subfloat[]{
      \includegraphics[width=7.5cm]{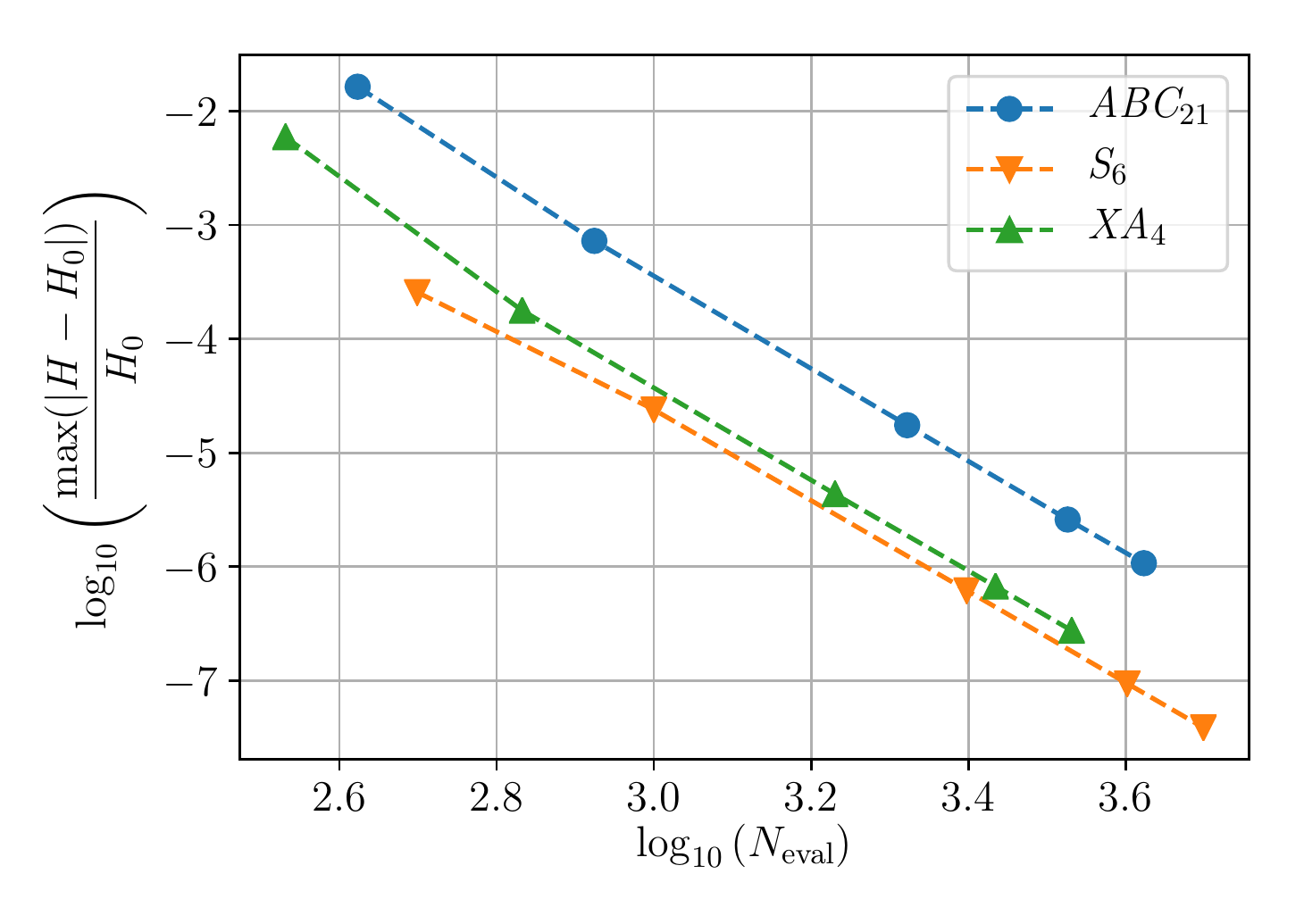}}
    \subfloat[]{
      \includegraphics[width=7.5cm]{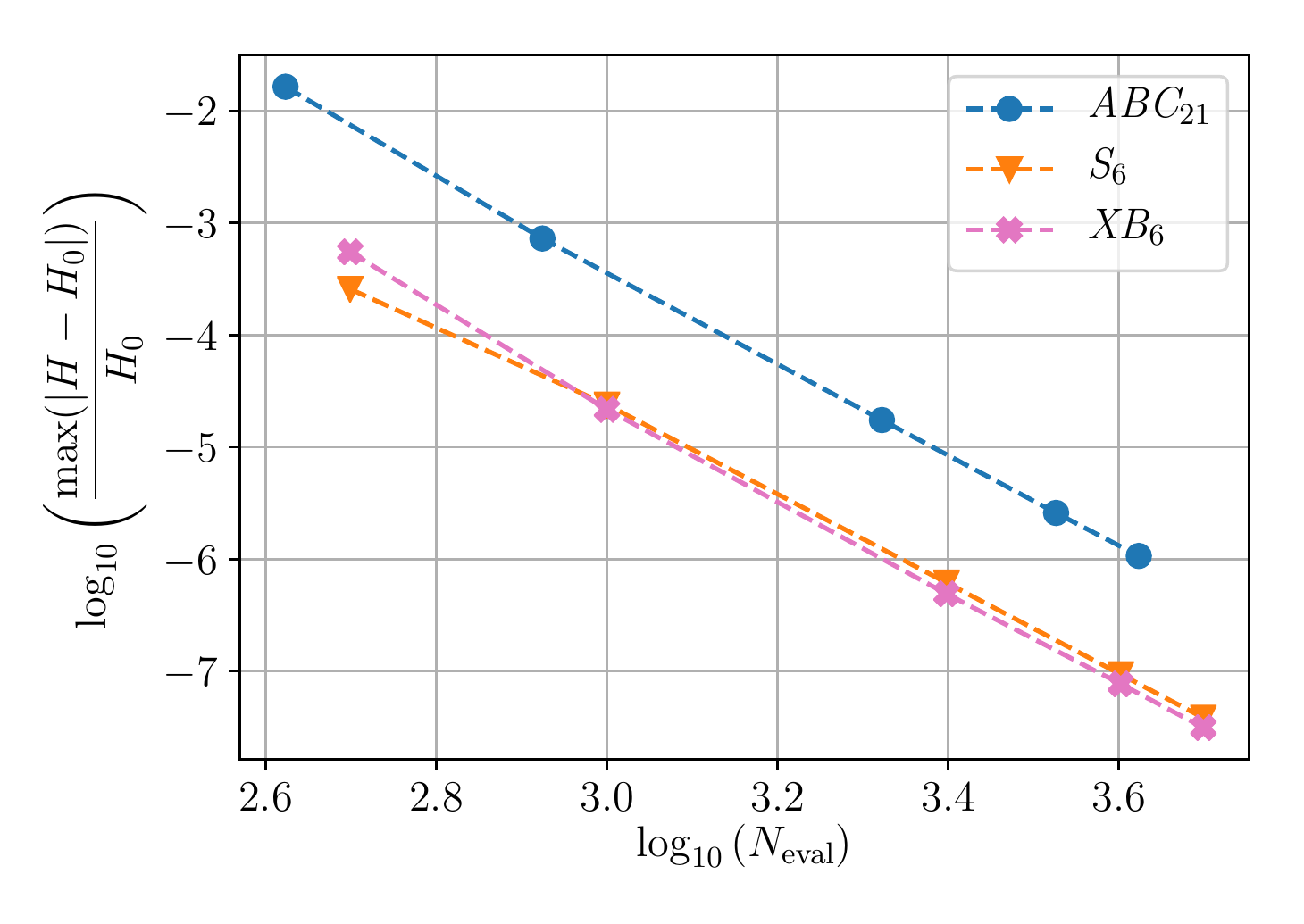}}
    \\ 
    \subfloat[]{
      \includegraphics[width=7.5cm]{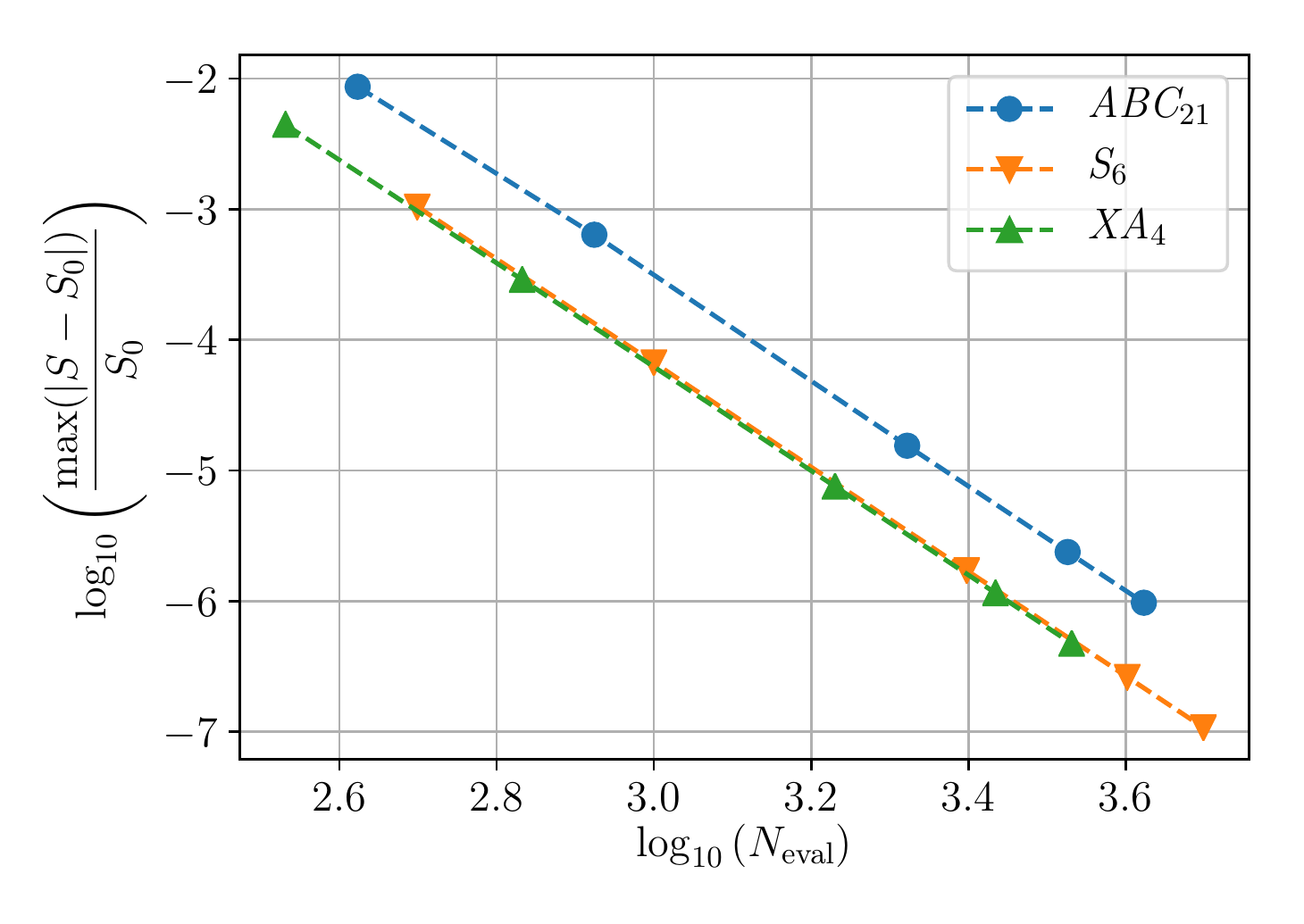}}
    \subfloat[]{
      \includegraphics[width=7.5cm]{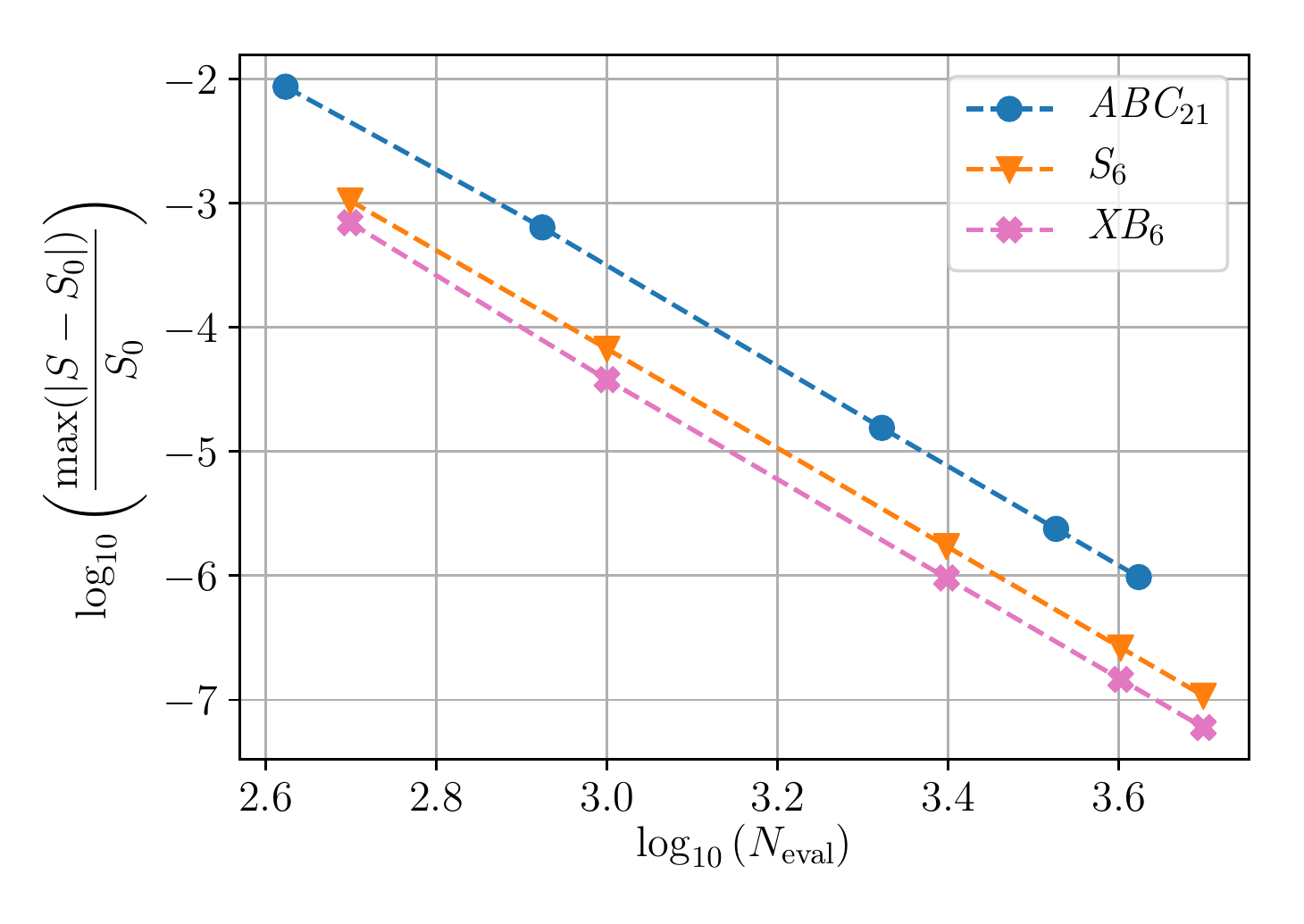}}
    \caption{\small Relative error in conserved quantities for the DDNLS system due to each of the best numerical methods tested. (\textbf{a}) Relative error in energy for \textit{XA}$_4$ compared to \textit{ABC}$_{21}$ and  \textit{S}$_6$. (\textbf{b}) Relative error in energy for \textit{XB}$_6$ compared to \textit{ABC}$_{21}$ and  \textit{S}$_6$. (\textbf{c}) Relative error in norm for \textit{XA}$_4$ compared to \textit{ABC}$_{21}$ and  \textit{S}$_6$. (\textbf{d}) Relative error in norm for \textit{XB}$_6$ compared to \textit{ABC}$_{21}$ and  \textit{S}$_6$.
      \label{fig:DDNLS1}}
  \end{center}
\end{figure}

\begin{figure}[!h]
  \begin{center}
    \subfloat[]{
      \includegraphics[width=7.5cm]{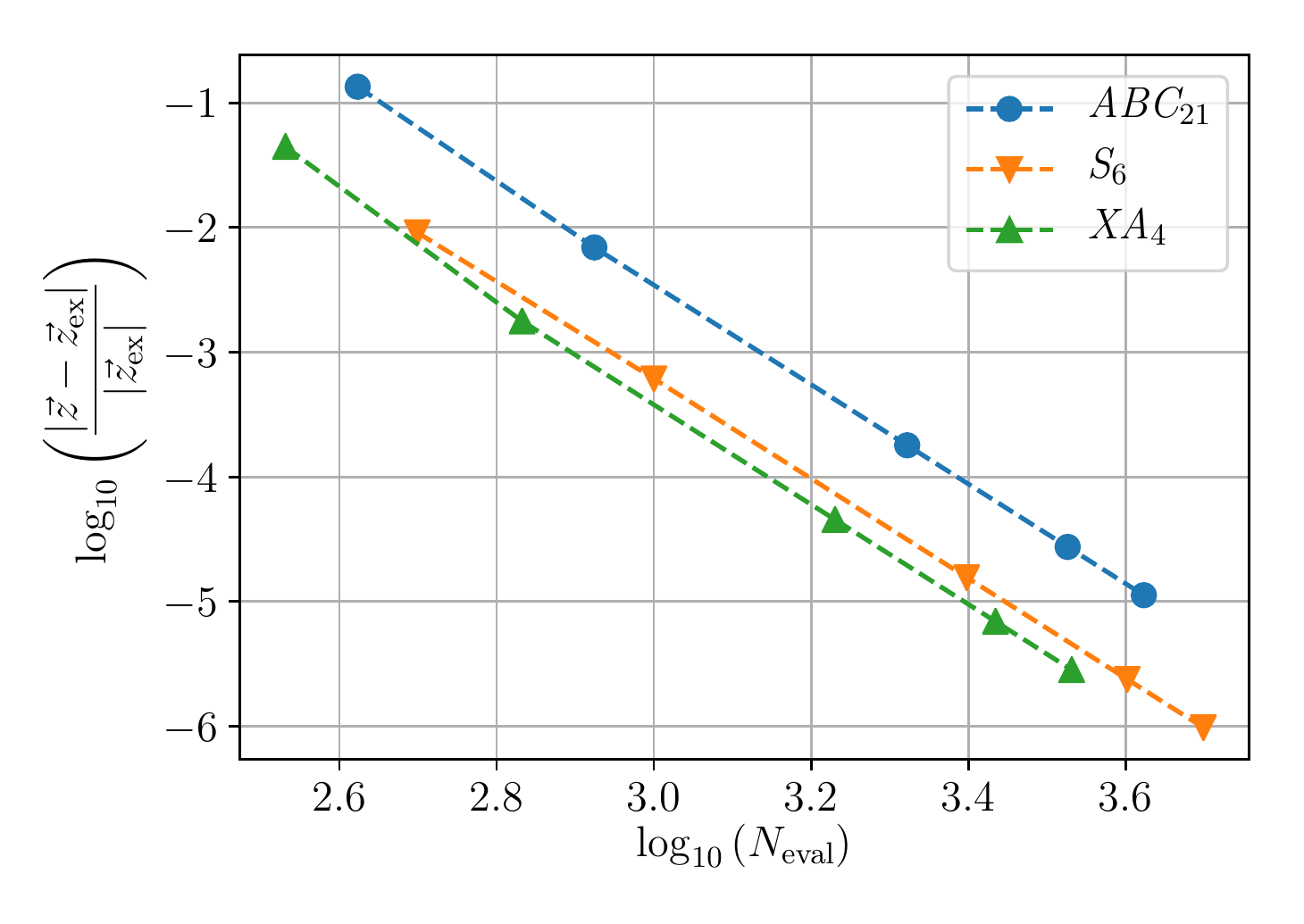}}
    \subfloat[]{
      \includegraphics[width=7.5cm]{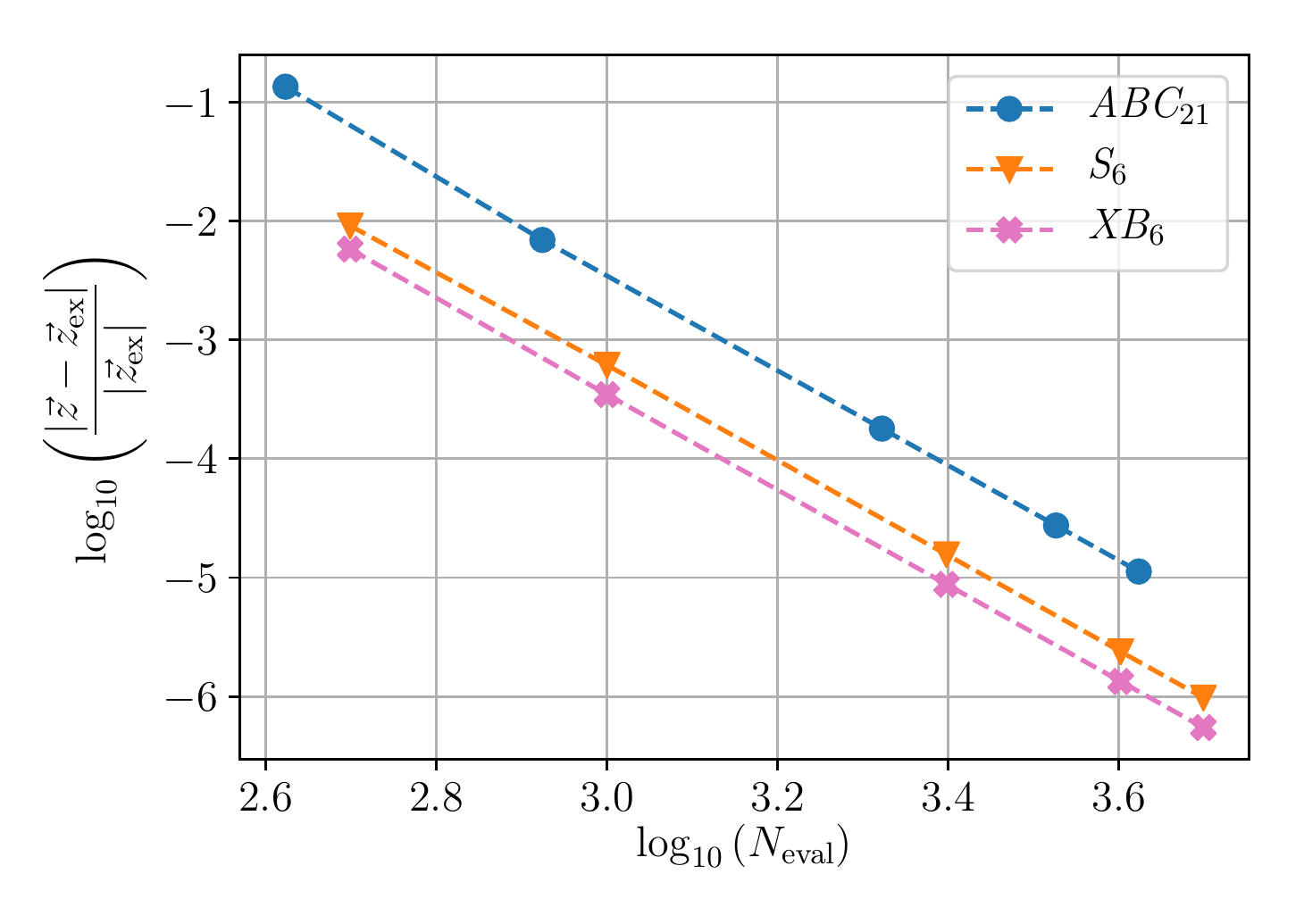}}
\caption{\small Relative error of trajectories for the DDNLS system. Same as in Figure \ref{fig:DDNLS1}. (\textbf{a}) Relative error of trajectories for \textit{XA}$_4$ compared to \textit{ABC}$_{21}$ and  \textit{S}$_6$. (\textbf{b}) Relative error of trajectories for \textit{XB}$_6$ compared to \textit{ABC}$_{21}$ and  \textit{S}$_6$.
\label{fig:DDNLS3}}
\end{center}
\end{figure}

To compare the performance of the numerical integrators previously considered, we take a lattice of $N= 1000$ sites and
fixed boundary conditions, $q_0 = p_0 = q_{N+1} = p_{N+1} = 0$. As in \cite{skokos14hot,danieli19ceo}, 
we~excite, at the initial time $t=0$, 21 central sites by taking the $q_i$ at 
random in the interval $[0,1]$ and the respective $p_i$ in such a way that each site has the same
constant norm 1, so that the total norm of the system is $S = 21$. Moreover, $\beta = 0.72$, $W=4$ and the random disorder parameters $\epsilon_j$ are chosen so that the total energy is $H \approx -29.63$. As in the previous example, we integrate until the final time $t_f = 10$ and compute the maximum relative error in energy and in norm along the integration interval. The results are depicted in Figure \ref{fig:DDNLS1}, with the top diagrams corresponding to the error in energy and the bottom to the error in norm. The same notation has been used for the tested methods. Finally, in Figure  \ref{fig:DDNLS3} we collect the error in the phase space. As before,
the reference solution is obtained with the DOP853 routine. Notice that for this non trivial example the new schemes \textit{XA}$_4$ and 
especially \textit{XB}$_6$ show a better efficiency than \textit{S}$_6$, and not only with respect to the preservation of the invariants, but also in the computation of trajectories.


\section{Concluding Remarks}
\label{sec.6}

In this work we have presented two different families of fourth-order composition methods especially designed for problems that can be separated into 
three parts in such a way that each part is explicitly solvable. In addition to the usual optimization criteria applied in the literature to choose the free parameters
in the composition, we have introduced another one especially oriented to problems where the energy is a constant of motion. The schemes constructed in this way
show an improved behavior, and in fact one of the methods exhibits a superior performance to the familiar scheme $S_6$ of Table \ref{table.2} on the 
tested examples. Other relevant examples include certain nonlinear oscillators, Poisson--Maxwell equations arising in plasma physics, and the
treatment of non separable Hamiltonian dynamical systems \cite{tao16esa}.

Although only problems separable into three parts have been considered here, it is clear that the schemes we have introduced can also be applied to differential
equations split into any number of pieces $n \ge 3$. The only modification one requires is to formulate the corresponding first order scheme $\chi_h$ and 
its adjoint $\chi_h^*$. One should be aware, however, that augmenting the number $n$ leads to evaluating an increasingly large number of flows for
methods with large values of $s$, with the subsequent deterioration in performance.

{
An important topic not addressed in this study concerns the stability of the proposed methods. Typically, for a  given method there exists a critical step size
$h^*$ such that it will be unstable for $|h| > h^*$. Of course, one is interested in methods with $h^*$ as large as possible. The linear stability of splitting
methods has been analyzed in particular in \cite{mclachlan97osp,blanes08otl}, where highly efficient schemes with optimal stability polynomials have 
presented for numerically approximating the evolution of linear problems. In the nonlinear case, however, the situation is more involved. In  \cite{mclachlan02foh},
a crude measure of the nonlinear stability of a given time symmetric scheme of order $r$ is proposed, taking into account the error terms of orders
$r+1$ and $r+3$. The stability of splitting methods in the particular setting of (semidiscretized) partial differential equations with stiff terms have been considered, 
in particular, in
\cite{hundsdorfer03nso,ropp05soo}. A theorem is presented \cite{ropp05soo} concerning
the stability of operator-splitting methods applied to linear reaction-diffusion equations with indefinite reaction terms which controls both low and high
wave number instabilities. In any case, this result only affects methods up to order 2 with real and positive coefficients, whereas the application of
splitting and composition methods of higher order with real coefficients in this setting leads to severe instabilities due to the existence of negative coefficients.
The methods we have presented here are aimed at non-stiff problems, and they do not exhibit, at least for the examples we have considered, special step
size restrictions in comparison with other splitting methods from the literature.}

Finally, it is worth remarking that the local error estimators for composition methods proposed in \cite{blanes19sac} based 
on the construction of lower order schemes obtained at each step as a linear combination of the intermediate stages of the main integrator, 
can also be used in this setting. As a consequence, it is quite
straightforward to implement the methods presented here with a variable step size strategy if necessary.



\vspace{6pt}


\subsection*{Acknowledgements}
This work has been funded by Ministerio de Econom\'{\i}a, Industria y Competitividad (Spain) through project MTM2016-77660-P (AEI/FE\-DER, UE)
and by Universitat Jaume I (projects UJI-B2019-17 and GACUJI/2020/05).
A.E.-T. has been additionally supported by the predoctoral contract BES-2017-079697 (Spain).

\end{document}